\documentclass[11pt,a4paper]{article}
\usepackage[utf8]{inputenc}
\usepackage[T1]{fontenc}
\usepackage{amsthm,mathrsfs,amsmath,amssymb}
\usepackage{dsfont}
\usepackage{comment}
\usepackage{indentfirst}
\usepackage{diagbox}
\usepackage{multirow}
\usepackage{booktabs}
\usepackage{cite}
\usepackage[colorlinks=true]{hyperref}
\hypersetup{urlcolor=blue, citecolor=red}
\numberwithin{equation}{section}
\allowdisplaybreaks
\makeatother
\date{}

\newtheorem{theorem}{Theorem}[section]
\newtheorem{lemma}[theorem]{Lemma}

\theoremstyle{definition}
\newtheorem{definition}[theorem]{Definition}
\newtheorem{remark}[theorem]{Remark}

\begin{document}
\baselineskip 15pt \setcounter{page}{1}
\title{  Normalized solutions for Schr\"{o}dinger-Maxwell equation with Sobolev critical   and mass subcritical perturbation\footnote{Supported by
National Natural Science Foundation of China (No.~12371120 and No.~12661024); and the National Natural Science Foundation of China under the Tianyuan Fund for Mathematics (No.~12526512); and the Science and Technology Research Program of Chongqing Municipal Education Commission (No. KJQN202601320).}}
\author{{Jin-Cai Kang$^a$,\ \ Yong-Yong Li$^b$,\ \ Chun-Lei Tang$^c$\footnote{Corresponding author. E-mail address:   tangcl@swu.edu.cn (C.-L. Tang)}}\\
{\small \emph{$^a$Chongqing Key Laboratory of Group} \& \emph{Graph Theories and Applications,}}\\
{\small\emph{Chongqing University of Arts and Sciences, Chongqing,  {\rm402160}, China}}\\
{\small\emph{$^b$College of Mathematics and Statistics, Northwest Normal University, Lanzhou, {\rm730070}, China}}\\
{\small\emph{$^c$School of Mathematics and Statistics, Southwest University, Chongqing, {\rm400715}, China}} }

\date{}
\maketitle

\baselineskip 16pt
\begin{quote}
  {\bf Abstract:}
  In this paper, we study the Schr\"{o}dinger-Maxwell equation with critical growth and mass subcritical perturbation
  \begin{align*}
  \begin{cases}
  \displaystyle - \Delta u+\lambda u+  \left(|x|^{-1}\ast| u|^2\right) u
  =|u|^{4}u +|u|^{p-2}u\quad\quad \text{in} \
  \mathbb{ R}^3,\\
  \displaystyle \int_{\mathbb{R}^3}|u|^2dx=a^2,
  \end{cases}
 \end{align*}
  where $a>0$ is prescribed,  $3<p<\frac{10}{3}$, $\ast$ denotes the convolution and $\lambda\in \mathbb{R}$ appears as a Lagrange multiplier.
 When the nonlocal convolution term is absent, the above equation can be reduced to the Schr\"{o}dinger equation. In this case,
the previous works showed that it admits a local minimizer (ground-state) and a mountain-pass type solution. By contrast, we surprisingly find that the above equation has no solutions with negative energy provided  $a>0$ is sufficiently small.
 More precisely, we rigorously prove that the above equation possesses   a mountain-pass type solution with positive energy   for sufficiently small $a>0$, which reveals an essential difference in solution structure between the Schr\"{o}dinger-Maxwell  equation and the classical Schr\"{o}dinger equation in the parameter range $3<p<\frac{10}{3}$.  The main difficulty encountered in this work is the loss of compactness arising from the Sobolev critical term.
  Furthermore, we   explore the dynamical properties of the corresponding time-dependent problem and establish the strong instability of its standing wave solutions.  It is worth pointing out that this work refines and extends the existing literature by covering the unexplored parameter regime $p\in(3,\frac{10}{3})$.

\end{quote}


\section{Introduction and main results}
\setcounter{section}{1}

In this paper, we investigate the standing waves   for a class of  Schr\"{o}dinger-Maxwell equations involving combined power nonlinearities
  \begin{align}\label{kk1}
   i \varphi_t+\Delta \varphi-\left(|x|^{-1}\ast| \varphi|^2\right)\varphi+|\varphi|^{p-2}\varphi+|\varphi|^{4}\varphi=0,
  \end{align}
where $3< p<\frac{10}{3}$ , $\varphi(x, t):\mathbb{R}^3\times [0, T]\rightarrow \mathbb{C} $ is the wave function   and  $(|x|^{-1}\ast  |\varphi| ^2)$ is a repulsive nonlocal Coulomb potential.  A class of equations with repulsive nonlocal Coulomb potentials can be derived via approximating the Hartree-Fock equation that describes multi-particle quantum mechanical systems (see \cite{2002RMP-BF, 1984CMP-L}). Specifically, the wave equation \eqref{kk1} admits combined nonlinearities consisting of several distinct nonlinear components, rather than a single nonlinear term. This composite nonlinear structure governs the intricate self-interaction and dynamical properties of the wave solutions. Tao and co-workers first conducted pioneering research on nonlinear Schr\"{o}dinger equations with combined nonlinearities  in \cite{2007CPDETao}. Their influential study has attracted considerable research attention and inspired numerous subsequent investigations in this area, see for example \cite{2017ARMA,  2013CMP} and the references therein.

The standing wave solution $\varphi(t, x) =e^{i \lambda t}u(x)$ of Eq.~\eqref{kk1} corresponds to a solution of the following   stationary elliptic equation
 \begin{align}\label{kkkk1}
 - \Delta u+\lambda u+ \left(|x|^{-1}\ast |u|^2\right)u
 =|u|^{p-2}u + |u|^{4}u\quad \quad \text{in} \ \mathbb{ R}^3.
 \end{align}
 For the resulting focusing variational problem \eqref{kkkk1}, existing studies mainly adopt two common research strategies. One is the fixed-frequency setting, where the frequency parameter $\lambda$ is fixed.   The other is the fixed-mass setting, where the $L^2$-norm of the solution is fixed. To date, most existing works focus on the fixed-frequency case, and comprehensive progress can be found in \cite{2010AIHP-ADP, 2006JFA,   2023-JGA-KLT, 2016-AMPA-LWZ} and the references therein. In comparison, the prescribed mass   problem of \eqref{kkkk1} remains less explored. Nevertheless, such prescribed mass   problems are physically meaningful, as they closely characterize the Bose-Einstein condensation process in quantum systems. This practical physical background strongly motivates us to study  Eq. \eqref{kkkk1} with a prescribed $L^2$-norm. In this framework, $\lambda$ is no longer a preset constant but serves as an unknown Lagrange multiplier.

In what follows, we introduce a generalized version of Schr\"{o}dinger-Maxwell equations   with combined power nonlinearities
 \begin{align}\label{L3}
  \begin{cases}
  \displaystyle - \Delta u+\lambda u+  \kappa\left(|x|^{-1}\ast|u|^2\right)u
  =|u|^{q-2}u +\mu|u|^{p-2}u\quad\quad \text{in} \
  \mathbb{R}^3,\\
  \displaystyle \int_{\mathbb{R}^3}|u|^2dx=a^2,
  \end{cases}
 \end{align}
where  $a, \kappa, \mu>0$, $2<p, q\leq6$   and $\lambda\in \mathbb{R}$ appears as a Lagrange multiplier.
To search for normalized solutions of Eq.~\eqref{L3}, we consider the critical points of the energy functional
\begin{align*}
 E(u)=\frac{1}{2}\int_{\mathbb{R}^3}|\nabla u|^2dx+\frac{\kappa}{4} \int_{\mathbb{R}^3} \int_{\mathbb{R}^3} \frac{|u(y)|^2} {|x-y|}dy|u(x)|^{2}dx
-\frac{1}{q}\int_{\mathbb{R}^3}|u|^{q}dx-\frac{1}{p}\mu\int_{\mathbb{R}^3}|u|^{p}dx
\end{align*}
on the constraint
$$
S_a:=\left\{u\in H^1(\mathbb{R}^3, \mathbb{C}): \int_{\mathbb{R}^3}|u|^2dx=a^2 \right\},
$$
To facilitate the subsequent discussion, we let $\overline{q}:=\frac{10}{3}$ and
 $$e(a):=\inf_{S_a} E(u).$$  When $\kappa=0$, Eq.~\eqref{L3} reduces to the Schr\"{o}dinger equation, which has   been extensively studied, see for example \cite{ 2020Soave, 2021-JFA-Soave, 2022JFAWei,2022.9}  and the references therein.
If $\kappa>0, \mu=0$,  $q\in(2, \overline{q})$, $E$ is bounded from below and coercive on $S_a$, and the existence of minimizer for $E$ constrained on $S_a$ has been studied in \cite{2011ZAMP-BS, 2013ZAMP-JL,   2011JFA-BS, 2004JSP}.   If $\kappa>0, \mu=0$, $q\in(\overline{q}, 6)$,  $E$ is unbounded from below on $S_a$.  Having established a bounded Palais-Smale sequence by using a mountain-pass argument on $S_a$, Bellazzini, Jeanjean and Luo in  \cite{2013PLMS} found critical points of the functional $E$ for $a> 0$ sufficiently small.  After then, if $\kappa>0, \mu=0$, $q=6$,    there are few results  presented in  \cite{2021-JDE-JL}. They showed that Eq.~\eqref{L3} does not admit positive solution in $H^1(\mathbb{R}^3)$. Moreover, for the general double power nonlinearity case where  $\kappa>0, \mu>0$,  $2<p<\frac{8}{3}< \overline{q}<q\leq 6$,  the previous work \cite{kang, 2025-JDE-Kang}  explored the existence of normalized solutions to Eq.~\eqref{L3}. Based on the Nehari-Poho\v{z}aev manifold technique, they successfully proved the existence of ground state solution and   mountain-pass type (shortly, M-P type) solution. Nevertheless, the case with $ \frac{8}{3}\leq p<\frac{10}{3}<q\leq6$ remains unresolved and lacks relevant research conclusions.  For more studies concerning normalized solutions of Schr\"{o}dinger-Maxwell equations, readers can refer to \cite{2017JMAA-ZZ, 2025-CV-P, 2014JMAA-L, 2018ZAMP-YL,  2020JMAA-CT, 2023IJM, 2026-Liu-JGA, 2026-Li-JGA} and the relevant references therein. To clearly demonstrate the research progress for Eq.~\eqref{L3}, we summarize all aforementioned existing results in Table \ref{T2}.

\begin{table}
\begin{small}
\begin{center}
\caption{Summary of results for Eq. \eqref{L3}}\label{T2}
\begin{tabular}{cccc}
\hline
   $\kappa, \mu, p, q$&Minimizer&  M-P type\\
\hline
  \multirow{2}{*}{$\kappa=0, \mu>0, 2<p<\overline{q}<q\leq6$}& \multirow{1}{*}{local minimizer (ground-state) exists}    & \multirow{1}{*}{M-P type}\\
  &\multirow{1}{*}{see \cite{  2020Soave, 2021-JFA-Soave, 2022JFAWei, 2022.9}   }&\multirow{1}{*}{see \cite{  2020Soave, 2021-JFA-Soave, 2022JFAWei, 2022.9}}\\
  \hline
   \multirow{2}{*}{$\kappa>0, \mu=0, 2<q<3$}& \multirow{1}{*}{if~$a>0$~small enough,   $e(a)<0$~is achieved}&\multirow{2}{*}{---}\\
     &\multirow{1}{*}{see \cite{2011JFA-BS}   }&\\
  \hline
    \multirow{2}{*}{$\kappa>0, \mu=0, q=\frac{8}{3}$}&  \multirow{1}{*}{if~$a>0$~suitably small,  $ e(a)<0$~is achieved}&\multirow{2}{*}{---}\\
     &\multirow{1}{*}{see \cite{2004JSP} }&\\
\hline
 \multirow{5}{*}{ $\kappa>0, \mu=0, q\in(3, \overline{q})$}&  \multirow{1}{*}{
if~$a\in(0,  \widetilde{a}),   e(a)=0$ is not achieved;}& \\
  &\multirow{1}{*}{if~$a=\widetilde{a},  e(a)=0$~is achieved;}&\multirow{2}{*}{---} \\
  &\multirow{1}{*}{ if~$a>\widetilde{a},  e(a)<0$~is achieved;}&\\
   &\multirow{1}{*}{if~$a>0$ is small,  $E|_{S_a}$ has no critical points}&\\
   &\multirow{1}{*}{see \cite{2011ZAMP-BS,  2013ZAMP-JL}}&\\
\hline
  \multirow{2}{*}{$\kappa>0, \mu=0, q=3$}&\multirow{1}{*}{$\forall a>0,  e(a)=0$~is not achieved};&  \\
   &\multirow{1}{*}{and $E|_{S_a}$ has no critical points, see \cite{2013ZAMP-JL} }&\\
\hline
  \multirow{4}{*}{$\kappa>0, \mu=0, q= \overline{q}$} & \multirow{1}{*}{if~$a\in(0, \overline{a}),  e(a)=0$~is not achieved;}&  \multirow{3}{*}{---}\\
 &\multirow{1}{*}{if~$a>\overline{a},  e(a)=-\infty$ is not achieved}& \\
 &\multirow{1}{*}{if~$a$ is small, $E|_{S_a}$ has no critical points}&\\
  &\multirow{1}{*}{see \cite{2013ZAMP-JL} }&\\
\hline
 \multirow{2}{*}{$\kappa>0, \mu=0, \overline{q}<q<6$}& \multirow{2}{*}{---} & \multirow{1}{*}{M-P type}  \\
  & & \multirow{1}{*}{see \cite{2013PLMS} }\\
\hline
  \multirow{2}{*}{$ \kappa>0, \mu=0, q=6$}& \multirow{1}{*}{$\forall a>0$,  $E|_{S_a}$ has no critical points } \\
  &   \multirow{1}{*}{see \cite{2021-JDE-JL} }\\
\hline
   \multirow{2}{*}{$\kappa>0, \mu>0, 2<p<\frac{8}{3}<\overline{q}<q\leq6$}& \multirow{1}{*}{local minimizer (ground-state) exists }   & \multirow{1}{*}{M-P type}\\
& \multirow{1}{*}{see \cite{kang, 2025-JDE-Kang} }& \multirow{1}{*}{see \cite{kang, 2025-JDE-Kang} }\\
\hline
\end{tabular}
\end{center}
\end{small}
\end{table}

Building on our previous research \cite{kang, 2025-JDE-Kang} and enlightened by  the research summary in Table \ref{T2}, this paper further explores the solvability of Eq.~\eqref{L3} by focusing on the uninvestigated parameter regime with $\kappa>0$, $\mu>0$, and $\frac{8}{3}<p< \overline{q}<q=6$. Beyond the existence of stationary normalized solutions, the dynamical stability and instability of solutions to the associated initial-value problem \eqref{kk1} also represent important unsolved issues in nonlinear analysis. Addressing these research gaps and continuing our prior investigations, this work is dedicated to studying the existence of normalized solutions for the following Schr\"{o}dinger-Maxwell equation  with combined nonlinearities
 \begin{align}\label{L1}\tag{${ \mathcal{P}}$}
  \begin{cases}
  \displaystyle - \Delta u+\lambda u+  \left(|x|^{-1}\ast| u|^2\right) u
  =|u|^{4}u +|u|^{p-2}u\quad \quad \text{in} \
  \mathbb{R}^3,\\
  \displaystyle \int_{\mathbb{R}^3}|u|^2dx=a^2,
  \end{cases}
 \end{align}
where  $a >0$,     $3<p<\frac{10}{3}$, $\ast$ denotes the convolution in $\mathbb{R}^3$ and $\lambda\in \mathbb{R}$ appears as an unknown  Lagrange multiplier.

For the convenience of narration, we let $\mathcal{H} := H^1(\mathbb{R}^3, \mathbb{C})$ denote the usual Sobolev space equipped with the norm
$$
\|u\| = \left( \int_{\mathbb{R}^3} |\nabla u|^2 + |u|^2 \, dx \right)^{\frac{1}{2}}.
$$
In addition, we define   $\mathcal{H}_r=\{ u\in\mathcal{H}: u(x)=u(|x |)\}$ and   $S_{a,r}:= S_a\cap  \mathcal{H}_r$.
As we all know, solutions of Eq. \eqref{L1} correspond to critical points of the energy functional $\mathcal{J}: \mathcal{H}\rightarrow \mathbb{R}$  on the constraint $S_a$, where $\mathcal{J}$ is given by
 \begin{align}\label{R5}
 \mathcal{J}(u)
 =\frac{1}{2}\int_{\mathbb{R}^3}|\nabla u|^2dx+\frac{1}{4} \int_{\mathbb{R}^3} \int_{\mathbb{R}^3} \frac{|u(y)|^2} {|x-y|}dy|u(x)|^{2}dx
 -\frac{1}{6}\int_{\mathbb{R}^3}|u|^{6}dx-\frac{1}{p}\int_{\mathbb{R}^3}|u|^{p}dx
\end{align}
 Clearly, $\mathcal{J}$  is well defined and of class $C^1$.

We will be particularly interested in ground state solutions, defined as
follows:
  \begin{definition}\label{de1.1}
  We say that  $\widetilde{u}$ is a ground state radial solution  of  Eq. \eqref{L1}  on $S_{a, r}$ if it is a solution to Eq. \eqref{L1} having minimal energy among all the solutions which belongs to $S_{a, r}$:
\begin{align*}
  d\mathcal{J}|_{S_{a, r}}(\widetilde{u})=0 \quad \text{and} \quad \mathcal{J}(\widetilde{u})=\inf\left\{\mathcal{J}(u):   d\mathcal{J}|_{S_{a, r}}(u)=0 \ \text{and}\ u\in S_{a, r} \right\}.
\end{align*}
The set of the radial ground states will be denoted by $Z_{a,r}$.
 \end{definition}

We also recall the notion of   instability  that we are interested in:

  \begin{definition}\label{de1.2}
 A standing wave $e^{i\lambda t}u$ is strongly unstable if for every $\varepsilon> 0$, there exists $\varphi_0 \in \mathcal{H}$ such
that $\|u- \varphi_0\|<\varepsilon$,
    and $\varphi(t,\cdot)$ blows-up in finite time, where $\varphi(t,\cdot)$ denotes the solution to \eqref{kk1} with initial datum  $\varphi_0$.
 \end{definition}

 Now, the main results in this paper are stated as follows:

\begin{theorem}\label{K-TH1}
Let    $p\in(3, \frac{10}{3} )$ and
\begin{align*}
\overline{a}:=\min\left\{a_\ast, \overline{a_\ast},  \left(\frac{4}{C_1(p)}\right)^{\frac{3p-8}{8(p-3)}}, \     \left(\frac{1}{C_2(p)}\right)^{\frac{3p-8}{8(p-3)}},    \   \left(\frac{1}{3C_3(p)}\right)^{\frac{3p-8}{8(p-3)}}, \   3^{\frac{1}{4}}\left(\frac{2}{3C_4(p)}\right)^{\frac{p-2}{2(p-3)}}\mathcal{S}^{-\frac{3}{4}}\right\},
\end{align*}
where $a_{\ast}, \overline{a_\ast}, \mathcal{S}$ and $ C_i(p) (i=1,2,3,4)$ will be defined in the later.   If $a\in(0, \overline{a})$,
then Eq. \eqref{L1}  possesses  a  mountain-pass type  normalized solution $(\lambda, u)\in  \mathbb{R}^+\times   \mathcal{H}_r$,
  where $u$ is a positive real-valued  function in $\mathbb{R}^3$.

\end{theorem}

\begin{remark}
A key observation of the present paper is that, for $p\in(3, \frac{10}{3} )$, Eq. \eqref{L1} exhibits geometric structures substantially different from those of the classical Schr\"{o}dinger equation with  Sobolev critical exponent and mass subcritical perturbation studied by Soave \cite{2021-JFA-Soave}, Wei et al.\cite{2022JFAWei} and Chen et al. \cite{2022.9}; see Lemma  \ref{jQ-Lem2.2}. More precisely,   existing works  \cite{2021-JFA-Soave, 2022JFAWei} verified that the Schr\"{o}dinger  equation with  Sobolev critical exponent and mass subcritical perturbation admits a ground state solution  with negative energy and a mountain-pass solution  with positive energy. In contrast, we show in this paper that any normalized solution to Eq. \eqref{L1}   must have positive energy.
As a byproduct,   Eq. \eqref{L1} with $p\in(3, \frac{10}{3} )$ has no normalized solutions with negative energy when $a\in (0, \overline{a})$.  In fact, the conclusion follows readily from \eqref{k5}.
\end{remark}

\begin{remark}
To prove Theorem \ref{K-TH1}, the main difficulty lies in verifying the compactness of Palais-Smale  sequences (shortly PS sequences), which arises from the Sobolev critical term. To resolve this obstacle, we first construct a special PS sequences that approximately satisfy Poho\v{z}aev identity  and thereby establish the boundedness of PS sequences. Next,  although the embedding $\mathcal{H}_r \hookrightarrow L^q(\mathbb{R}^3)$ with $q\in(2, 6)$ is compact, we still cannot derive the compactness of PS sequences in $L^2(\mathbb{R}^3)$. To overcome it, we introduce two key steps to establish compactness: (1) Owing to the Sobolev critical exponent, we derive an upper bound estimate for the mountain-pass level. Notice that, this threshold differs from those adopted in existing literature (see Lemma \ref{R-lem2}). (2)  We need to show that $\lambda>0$. It is challenging to obtain this result via the estimates and methods adopted in previous works. Therefore, in this paper, we utilize refined analytical techniques and establish several new estimates to achieve this goal (see Lemma \ref{R-lem6}).
\end{remark}

\begin{remark}
It is worth noting that the present paper constitutes a further investigation of our earlier work \cite{kang, 2025-JDE-Kang}. In that work, we proved that Eq. \eqref{L1} with $p\in(2,\frac{8}{3})$ admits a negative-energy ground state solution and a positive-energy mountain pass type solution. In particular, the question of whether Eq. \eqref{L1} possesses normalized solutions for the parameter range $p\in[\frac{8}{3}, 3]$ remains open.
\end{remark}

Next, we further study the dynamical behavior of solutions to problem \eqref{kk1}, with a particular focus on the strong instability of standing waves.
According to \cite{2003-C},   we know that the  initial-value problem \eqref{kk1} with initial datum   is  local well-posed   on $(-T_{min}, T_{max})$ with $T_{min}, T_{max}>0 $. Let $ \phi(t,x)$   be a solution of  the  initial-value problem \eqref{kk1}
with initial datum $u_0$ on $(-T_{min}, T_{max})$, which satisfies
\begin{align*}
|u_0|_2=| \phi|_2\quad\quad \text{and} \quad \quad \mathcal{J}(u_0)=\mathcal{J}( \phi).
\end{align*}
Moreover, by the rotational invariance of  problem \eqref{kk1},   similar to \cite[Theorem 4]{2013APDE},   we know that if the initial datum $u_0$ is radial, then $ \phi(t,x)$ is radial for all times.


\begin{theorem}\label{TH4}
Under the assumptions of Theorem \ref{K-TH1},
   let $u_0 \in S_{a, r}$ be such that $\mathcal{J}(u_0)<\inf_{u\in \mathcal{P}_{a}}\mathcal{J}(u)$, and if  $|x|u_0 \in L^2(\mathbb{R}^3, \mathbb{C})$ and $s_{u_0}<0$, where  $\mathcal{P}_{a}$  is defined in \eqref{k3} and $s_{u_0}$ is given in Lemma \ref{Q-Lem2.3}, then the solution $\phi$ of problem \eqref{kk1}  with initial datum $u_0$ blows-up
in finite time.
\end{theorem}

\begin{theorem}\label{TH5}
Under the assumptions of Theorem \ref{K-TH1},
  if $\widehat{u}$ is a ground state radial solution, then the associated Lagrange multiplier $\widehat{\lambda}$ is positive, and the
standing wave $e^{i\widehat{\lambda} t}\widehat{u} $ is strongly unstable.
\end{theorem}
\begin{remark}
The proof frameworks of Theorems \ref{TH4} and \ref{TH5} are inspired by the work of Soave \cite{2020Soave}; we implement tailored technical arguments adapted to our equation in the detailed derivation.
\end{remark}

The organization of this paper is as follows. Sec.~\ref{sec2} gives some  preliminaries.
In Sec.~\ref{sec3}, we prove Theorem \ref{K-TH1}.  Sec.~\ref{sec4} is devoted to the proof of Theorems \ref{TH4} and \ref{TH5}.

We conclude this section with the following notations applied throughout this paper:
\begin{itemize}
  \setlength{\itemsep}{-0.5pt}
  \item $\mathcal{D}^{1,2}(\mathbb{R}^3):= \mathcal{D}^{1,2}(\mathbb{R}^3,\mathbb{C})$ is the usual Sobolev space with the norm
  $\|u\|_{\mathcal{D}^{1,2}  }=\left(\int_{\mathbb{R}^3} |\nabla u|^2 dx\right)^{\frac12}$.
  \item $L^q(\mathbb{R}^3)= L^q(\mathbb{R}^3, \mathbb{C})$ is the Lebesgue space with the norm $|u|_q=\left(\int_{\mathbb{R}^3}|u|^qdx\right)^{\frac{1}{q}}$ for $q\in[1,+\infty)$.
      \item The sharp constant $\mathcal{S} =\inf \left\{\|u\|_{\mathcal{D}^{1,2}  }^2:u\in \mathcal{D}^{1,2}(\mathbb{R}^3)\ \mbox{and}\ |u|_6=1\right\}$.
  \item 
   $\mathbb{R}^+=(0, +\infty)$.
  \item $C$, $C_i,  i=1,2, \ldots $,   denote positive constant which may depend on  $p$, and may possibly various in different places.
  \item $B(y, r):=\left\{x\in \mathbb{ R}^3:|x-y|<r\right\}$ for $y\in\mathbb{R}^3$ and $r>0$.
\end{itemize}

\section{Preliminaries}\label{sec2}
In this section, we first recall the Gagliardo-Nirenberg inequality and several key properties of the nonlocal convolution term, which are essential to the analysis of this paper.

\begin{lemma}(see \cite{1983W})
 For any $u\in  \mathcal{H}$,   the following Gagliardo-Nirenberg inequality  holds
 \begin{align}\label{k10}
|u|_t\leq C_{t}|u|_2^{1-\gamma_t}|\nabla u|_2^{\gamma_t},
 \end{align}
  where $2< t <6$, the sharp constant $C_{t}$ is given by
 \begin{align*}
C_{t}^t=\frac{2t}{6-t}\left[\frac{6-t}{3(t-2)}  \right]^{\frac{3(t-2)}{4}}\frac{1}{|Q_t|_2^{t-2}},
 \end{align*}
and $Q_t$ is the unique positive radial solution of the equation $-\Delta Q+Q=|Q|^{t-2}Q$ in ${\mathbb{R}}^3$.
\end{lemma}

According to \cite{2006JFA}, one can present the following properties for the nonlocal term.
\begin{lemma}\label{gle4}
The following results hold:
\begin{description}
  \item [$(1)$]   $\displaystyle \int_{\mathbb{R}^3} \frac{|u(y)|^2} {|x-y|}dy\ge0$  for any $x\in{\mathbb{R}^{3}}$;
  \item [$(2)$] there exist $C_1,C_2>0$ such that $\displaystyle\int_{\mathbb{R}^{3}}\int_{\mathbb{R}^3} \frac{|u(y)|^2} {|x-y|}dy |u|^2dx\le{C_1|u|^4_\frac{12}{5}}\le{C_2\|u\|^4};$
  \item [$(3)$] if $u_n\rightarrow u$ in $L^{\frac{12}{5}}(\mathbb{R}^3)$, then $\displaystyle\int_{\mathbb{R}^3} \frac{|u_n(y)|^2} {|x-y|}dy\rightarrow \int_{\mathbb{R}^3} \frac{|u(y)|^2} {|x-y|}dy$ in $\mathcal{D}^{1,2}(\mathbb{R}^{3});$
  \item [$(4)$]    if $u_n\rightharpoonup u$ in $ \mathcal{H}$, then $\displaystyle\int_{\mathbb{R}^3} \frac{|u_n(y)|^2} {|x-y|}dy\rightharpoonup \int_{\mathbb{R}^3} \frac{|u(y)|^2} {|x-y|}dy$ in $\mathcal{D}^{1,2}(\mathbb{R}^{3})$.
\end{description}
\end{lemma}

Similar to \cite[Lemma 2.2]{2013ZAMP-JL} and \cite[Lemma 2.1]{2013-3M}, one obtains the following lemma.
\begin{lemma}\label{Q-Lem2.1}
When $p \in [3 ,4)$, there exists a constant $\overline{C}>0$, depending only on $p$, such that, for any
$u\in S_a$,
\begin{align*}
\int_{\mathbb{R}^3} \int_{\mathbb{R}^3} \frac{ |u(y)|^2|u(x)|^2}{|x-y|}  dy dx
\geq   \overline{C} \frac{|u|_p^{\frac{2p}{4-p}}}{ |\nabla u|_2^{\frac{2(2p-5)}{4-p}} |u|_2^{\frac{2(p-3)}{4-p}}}.
\end{align*}
\end{lemma}

\begin{proof}
Given $p\in[3,4)$, we apply interpolation inequalities to derive
\begin{align}\label{Q20}
|u|_p^p \le |u|_3^{3(4-p)}|u|_4^{4(p-3)}.
\end{align}
Furthermore, since $\bigl(|x|^{-1}*|u|^2\bigr)\in \mathcal{D}^{1,2}(\mathbb R^3)$ solves the equation
\begin{align}\label{Q21}
-\Delta \Phi = 4\pi|u|^2 \quad\quad\text{in }\mathbb R^3,
\end{align}
 multiplying \eqref{Q21} by $\bigl(|x|^{-1}*|u|^2\bigr)\in \mathcal{D}^{1,2}(\mathbb R^3)$ and integrating, we get
\begin{align*}
  \int_{\mathbb R^3}\bigl(|x|^{-1}*|u|^2\bigr)|u|^2 dx
=\frac{1}{4\pi}\int_{\mathbb R^3}\bigl|\nabla\bigl(|x|^{-1}*|u|^2\bigr)\bigr|^2 dx.
\end{align*}
Besides, multiplying \eqref{Q21} by $|u|$ and integrating, one has, for any $\eta>0$,
\begin{align*}
 4\pi \eta \int_{\mathbb R^3}|u|^3 dx
&=\eta \int_{\mathbb R^3}-\Delta\bigl(|x|^{-1}*|u|^2\bigr)|u|dx \\
&\le \eta \int_{\mathbb R^3}\nabla\bigl(|x|^{-1}*|u|^2\bigr)\cdot \nabla|u|dx \\
&\le \frac{1}{2}\int_{\mathbb R^3}\bigl|\nabla\bigl(|x|^{-1}*|u|^2\bigr)\bigr|^2 dx
+\frac{\eta^2}{2}\int_{\mathbb R^3}|\nabla u|^2 dx.
\end{align*}
After optimizing on $\eta $, we obtain that
\begin{align}\label{Q23}
|u|_{3}^6 \le  |\nabla u|_{2}^2 \int_{\mathbb R^3}\bigl|\nabla\bigl(|x|^{-1}*|u|^2\bigr)\bigr|^2 dx.
\end{align}
 Now, using Gagliardo-Nirenberg  inequality \eqref{k10}, there exists a constant $C>0$, depending only on $p$, such that
\begin{align}\label{Q24}
\int_{\mathbb R^3}|u|^4 dx \le C|\nabla u|_2^3|u|_2.
\end{align}
Substituting \eqref{Q23} and \eqref{Q24} into \eqref{Q20} yields
 \begin{align*}
 |u|_p^p\leq C |\nabla u|_2^{2p-5} \left(\int_{\mathbb{R}^3} \int_{\mathbb{R}^3} \frac{ |u(y)|^2|u(x)|^2}{|x-y|}  dy dx\right)^{\frac{4-p}{2}}|u
|_2^{p-3},
 \end{align*}
which gives the desired result.
\end{proof}

For any $u\in S_a$, let $(s\star u)(\cdot)=e^{\frac{3s}{2}}u(e^s\cdot)$,
one has $s\star u\in S_a$ and
 \begin{align}\label{k8}
 \psi_u(s):= \mathcal{J}(s\star u )
 =&\frac{e^{2s}}{2}\int_{\mathbb{R}^3}|\nabla u|^2dx+\frac{ e^{s}}{4} \int_{\mathbb{R}^3} \int_{\mathbb{R}^3} \frac{|u(y)|^2|u(x)|^{2}} {|x-y|}dydx \nonumber\\
 & -\frac{e^{6s}}{6}\int_{\mathbb{R}^3}|u|^{6}dx
 -\frac{e^{p\gamma_ps}}{p}
 \int_{\mathbb{R}^3}|u|^{p}dx,
\end{align}
where $\gamma_p=\frac{3(p-2)}{2p}$. Note that $p \gamma_p\in(\frac{3}{2}, 2)$ if $p\in(3, \frac{10}{3})$.
If $u\in  \mathcal{H}$ is a solution of system \eqref{L1}, then  the following Nehari  identity holds
 \begin{align}\label{kk4}
\int_{\mathbb{R}^3}|\nabla u|^2+\lambda\int_{\mathbb{R}^3}  |u|^2dx+ \int_{\mathbb{R}^3} \int_{\mathbb{R}^3} \frac{|u(y)|^2|u(x)|^{2}} {|x-y|}dydx  - \int_{\mathbb{R}^3}|u|^{6}dx- \int_{\mathbb{R}^3}|u|^{p}dx=0.
\end{align}
By \cite[(4.35)]{2004-ANS-D} (see also \cite[Theorem 2.2]{2006JFA}),    if $u\in  \mathcal{H}$ is a solution of  system \eqref{L1},
then $u$ satisfies the following Poho\v{z}aev identity:
 \begin{align}\label{kkk4}
\frac{1}{2}\int_{\mathbb{R}^3}|\nabla u|^2dx +\frac{3\lambda }{2} \int_{\mathbb{R}^3}  |u|^2dx+\frac{5}{4}\int_{\mathbb{R}^3} \int_{\mathbb{R}^3} \frac{|u(y)|^2|u(x)|^{2}} {|x-y|}dydx=
\frac{1}{2}\int_{\mathbb{R}^3}|u|^{6}dx+\frac{ 3}{p} \int_{\mathbb{R}^3}|u|^{p}dx.
\end{align}
Then,   combining   \eqref{kk4} and \eqref{kkk4},   any solution $u$ of system \eqref{L1} satisfies
\begin{align}\label{k4}
P(u):=\int_{\mathbb{R}^3}|\nabla u|^2dx+\frac{1}{4}\int_{\mathbb{R}^3} \int_{\mathbb{R}^3} \frac{|u(y)|^2|u(x)|^{2}} {|x-y|}dydx
- \int_{\mathbb{R}^3}|u|^{6}dx
-\gamma_p\int_{\mathbb{R}^3}|u|^{p}dx=0.
 \end{align}
Usually, \eqref{k4} is widely used in studying prescribed mass problem and called Poho\v{z}aev identity.

  To find the normalized solutions of Eq. \eqref{L1}, we introduce the  Poho\v{z}aev constrained set
\begin{align}\label{k3}
\mathcal{P}_{a} = \bigl\{ u \in \mathcal{H}_r : P(u) = 0 \bigr\} \cap S_{a,r}.
\end{align}
Obviously, the set $\mathcal{P}_{a}$ contains all of radial normalized solutions for Eq. \eqref{L1}. Therefore, this paper mainly studies the following minimization problem
\[
m_{a} = \inf_{u \in \mathcal{P}_{a}} \mathcal{J}(u).
\]
By simple direct calculation, one has
 \begin{align}\label{j-Q30}
P(s\star u)=e^{2s}\int_{\mathbb{R}^3}|\nabla u|^2dx+\frac{e^{s}}{4} \int_{\mathbb{R}^3} \int_{\mathbb{R}^3} \frac{|u(y)|^2|u(x)|^{2}} {|x-y|}dydx
-  e^{6s} \int_{\mathbb{R}^3}|u|^{6}dx
-\gamma_p e^{p\gamma_ps}
\int_{\mathbb{R}^3}|u|^{p}dx.
\end{align}
Clearly, $\psi_u'(s) = P(s \star u)$, which implies that for any $u \in S_{a, r}$, $s \in \mathbb{R}$ is a critical point for $\psi_u(s)$ if and only if $s \star u \in \mathcal{P}_{a}$. To further explore the geometric structure of  $\mathcal{J}|_{\mathcal{P}_{a}}$,  we consider the decomposition of $\mathcal{P}_{a}$ into the disjoint union
$$
\mathcal{P}_{a} = (\mathcal{P}_{a})^+ \cup (\mathcal{P}_{a})^0 \cup (\mathcal{P}_{a})^-,
$$
where
$$
(\mathcal{P}_{a})^{+(\text{resp. }0,-)} = \bigl\{ u \in \mathcal{P}_{a} : \psi_u''(0) > (\text{resp. }=,<)0 \bigr\}.
$$

Now, we are ready to establish the following lemmas characterizing the component of $\mathcal{P}_a$.

\begin{lemma}\label{jQ-Lem2.2}
 Let   $p\in(3, \frac{10}{3} )$ and $a\in(0, \overline{a})$.  Then $\mathcal{P}_a  =  \mathcal{P}_a^-$.
\end{lemma}
\begin{proof}
 For any $u \in \mathcal{P}_a$, then $P(u) = 0$, namely,
\begin{align*}
P(u)=\int_{\mathbb{R}^3}|\nabla u|^2dx+\frac{1}{4}\int_{\mathbb{R}^3} \int_{\mathbb{R}^3} \frac{|u(y)|^2|u(x)|^{2}} {|x-y|}dydx
- \int_{\mathbb{R}^3}|u|^{6}dx
-\gamma_p\int_{\mathbb{R}^3}|u|^{p}dx=0.
 \end{align*}
where $\gamma_p = \dfrac{3(p-2)}{2p}$. From \eqref{k8}, one has
$$
\psi_u''(0) = 2 \int_{\mathbb{R}^3} |\nabla u|^2 dx + \frac{1}{4}\int_{\mathbb{R}^3} \int_{\mathbb{R}^3} \frac{|u(y)|^2|u(x)|^{2}} {|x-y|}dydx - 6 |u|_6^6 -  p\gamma_p^2 |u|_p^p.
$$
We recall from \cite[Lemma 2.1]{2013-3M} that, for any $p\in(3, \frac{10}{3})$ and $u\in \mathcal{H}$,
\begin{align}\label{j-Q46}
\int_{\mathbb{R}^3} \int_{\mathbb{R}^3} \frac{|u(y)|^2|u(x)|^{2}} {|x-y|}dydx\geq
\widehat{C} \frac{|u|_p^{\frac{2p}{10-3p}}}{ |\nabla u|_2^{\frac{2(3p-8)}{10-3p}} |u|_2^{\frac{8(p-3)}{10-3p}}}.
\end{align}
Then, by \eqref{j-Q46} we deduce
\begin{align}\label{Q60}
\psi_u''(0) &= \psi_u''(0) - 6 P(u) \nonumber\\
&=-4 |\nabla u|_2^2- \frac{5}{4} \int_{\mathbb{R}^3} \int_{\mathbb{R}^3} \frac{|u(y)|^{2}|u(x)|^2} {|x-y|}dydx + \gamma_p(6-p\gamma_p) |u|_p^{p}  \nonumber\\
&\leq -4 |\nabla u|_2^2- \frac{5}{4} \widehat{C} \frac{|u|_p^{\frac{2p}{10-3p}}}{ |\nabla u|_2^{\frac{2(3p-8)}{10-3p}} |u|_2^{\frac{8(p-3)}{10-3p}}}+ \gamma_p(6-p\gamma_p) |u|_p^{p}
\end{align}
Let $\mathcal{A}:= \frac{5}{4} \widehat{C}\Big(|\nabla u|_2^{\frac{2(3p-8)}{10-3p}} |u|_2^{\frac{8(p-3)}{10-3p}}\Big)^{-1} $.
 Consider the auxiliary function
 \begin{align*}
f(x)=-4|\nabla u|_2^2 - \mathcal{A} x^{\frac{2}{10-3p}}+\gamma_p(6-p\gamma_p)x, \quad\quad x>0.
 \end{align*}
  By a simple calculation, one has
  \begin{align*}
f'(x) &= -\frac{2}{10-3p}\mathcal{A} \cdot x^{\frac{3p-8}{10-3p}} + \gamma_p(6-p\gamma_p),\\
f''(x) &= -\frac{2}{10-3p} \cdot  \frac{3p-8}{10-3p}  \mathcal{A} \cdot x^{\frac{6(p-3)}{10-3p}}< 0, \quad \quad \text{for all } x>0.
\end{align*}
Therefore, the function $f(x)$ has the unique global maximum at
\begin{align*}
\bar{x}=\left(\frac{\gamma_p(6-p\gamma_p)(10-3p)}{2\mathcal{A}} \right)^{\frac{10-3p}{3p-8}},
  \end{align*}
  and then
  \begin{align*}
  \max_{x>0}f(x)
  =f(\bar{x})
  =& -4|\nabla u|_2^2 - \mathcal{A} \overline{x}^{\frac{2}{10-3p}}+\gamma_p(6-p\gamma_p)\overline{x}\nonumber\\
=& -4|\nabla u|_2^2  - \mathcal{A}   \left( \frac{\gamma_p(6-p\gamma_p)(10-3p)}{2\mathcal{A}} \right)^{\frac{2}{3p-8}}\nonumber\\
&+\gamma_p(6-p\gamma_p)\left( \frac{\gamma_p(6-p\gamma_p)(10-3p)}{2\mathcal{A}} \right)^{\frac{10-3p}{3p-8}} \nonumber\\
= &-4|\nabla u|_2^2+\frac{3p-8}{2}\left(\frac{1}{\mathcal{A}}\right)^{\frac{10-3p}{3p-8}}\left( \gamma_p(6-p\gamma_p)\right)^{\frac{2}{3p-8}}\left(\frac{10-3p}{2}\right)^{\frac{10-3p}{3p-8}}\nonumber\\
=&  -4|\nabla u|_2^2 + C_1(p) |\nabla u|_2^2|u|_2^{\frac{8(p-3)}{3p-8}},
    \end{align*}
    where
    \begin{align}\label{k7}
    C_1(p):=\frac{3p-8}{2}\left[\gamma_p(6-p\gamma_p)\right]^{\frac{2}{3p-8}}\left(\frac{2(10-3p)}{5\widehat{C}}\right)^{\frac{10-3p}{3p-8}} \end{align}
  Back to \eqref{Q60},   one gets
 \begin{align}
\psi_u''(0) &\leq -4 |\nabla u|_2^2- \frac{5}{4} \widehat{C} \frac{|u|_p^{\frac{2p}{10-3p}}}{ |\nabla u|_2^{\frac{2(3p-8)}{10-3p}} |u|_2^{\frac{8(p-3)}{10-3p}}}+ \gamma_p(6-p\gamma_p) |u|_p^{p}\nonumber\\
&\leq  -4|\nabla u|_2^2  + C(p) |\nabla u|_2^2|u|_2^{\frac{8(p-3)}{3p-8}}\nonumber\\
 &\leq  \left(-4+C_1(p)a^{\frac{8(p-3)}{3p-8}} \right) |\nabla u|_2^2\label{Q61}\\
 &<0\nonumber
\end{align}
when $a<\left(\frac{4}{C_1(p)}\right)^{\frac{3p-8}{8(p-3)}}$.
This yields that $u \in \mathcal{P}_a^-$. Thus,  $\mathcal{P}_a = \mathcal{P}_a^-$ holds. We complete the proof.
\end{proof}

\begin{lemma}\label{Q-Lem2.3}
 Let   $p\in(3, \frac{10}{3} )$ and $a\in(0, \overline{a})$. For any $u \in S_{a, r}$, the function $\psi_u$ has a unique critical point $s_u$. In other words, there exists a unique $s_u \in \mathbb{R}$ such that $s_u \star u \in \mathcal{P}_a$. Moreover, $\mathcal{J}(s_u \star u) = \max_{s\in \mathbb{R}} \mathcal{J}(s \star u)$ and $\psi_u$ is strictly decreasing and concave on $(s_u, +\infty)$. In particular, $s_u < 0$ if and only if $P(u) < 0$. And the map $u \in S_{a, r} \mapsto s_u \in \mathbb{R}$ is of class $C^1$.
\end{lemma}

 \begin{proof}
Since $ \psi_u'(s)= P(s\star u )$, it suffices to verify that $\psi_u'(s)$ possesses exactly one root over $\mathbb{R}$. From \eqref{j-Q30}, for all $u\in S_{a, r}$, we obtain
\begin{align} \label{h20}
\psi_u'(s)&=
e^{2s}\int_{\mathbb{R}^3}|\nabla u|^2dx+\frac{e^{s}}{4} \int_{\mathbb{R}^3} \int_{\mathbb{R}^3} \frac{|u(y)|^2|u(x)|^{2}} {|x-y|}dydx
-  e^{6s} \int_{\mathbb{R}^3}|u|^{6}dx
-\gamma_p e^{p\gamma_ps}
\int_{\mathbb{R}^3}|u|^{p}dx.
 \end{align}
 For $p\in(3, \frac{10}{3} )$, we conclude $ \psi_u'(s)>0$ as $s\rightarrow-\infty$.  Accordingly, one can find $s_0\in \mathbb{R}$ such that   $ \psi_u(s)$ is strictly increasing on $(-\infty, s_0 )$. Meanwhile, $\psi_u(s)\rightarrow -\infty $  as $s\rightarrow+\infty$.   Consequently,  there is $s_1\in \mathbb{R}$ with $s_1>s_0$ such that
 \begin{align}\label{h36}
  \psi_u(s_1)=\max_{s\in \mathbb{R}} \psi_u(s).
  \end{align}
  Therefore, $ \psi_u'(s_1)=0$, namely, $s_1\star u\in \mathcal{P}_a$. Assume that there exists $s_1\neq s_2\in \mathbb{R}$ satisfying   $s_2\star u\in \mathcal{P}_a$.   Without loss of generality, assume that $ s_2>s_1$.  By Lemma \ref{jQ-Lem2.2}, one infers $ \psi_u''(s_2)<0$. Hence, there is $s_3\in(s_1, s_2)$ satisfying
  $$
  \psi_u(s_3)=\min_{s\in(s_1, s_2) } \psi_u(s).
  $$
  Thus, we have $ \psi_u''(s_3)\geq0$ and $ \psi_u'(s_3)=0$, which implies that $s_3\star u\in \mathcal{P}_a$ and      $s_3\star u\in (\mathcal{P}_a)^+\cup (\mathcal{P}_a)^0$. This is a contradiction. Hence,  setting $s_u=s_1$, we know that $s_u\in \mathbb{R}$ is an unique number satisfying $ s_u\star u \in \mathcal{P}_a$. And then  it follows from \eqref{h36} that $\mathcal{J}(s_u\star u )=\max_{s\in \mathbb{R}}\mathcal{J}(s\star u ) $. Furthermore,
  \begin{align} \label{h47}
  \psi_u'(s)<0  \quad \Leftrightarrow \quad  s>s_u.
\end{align}
 So,   $ \psi_u'(0)=P(u)<0$ if and only if $s_u<0$.  Next, we show that   $ \psi_u$   is strictly decreasing and concave on  $(s_u, +\infty)$. Clearly, it follows from \eqref{h47} that  $ \psi_u$   is strictly decreasing  on  $(s_u, +\infty)$. Now, we   verify that $\psi_u''(s)<0$ on  $(s_u, +\infty)$.
   From \eqref{j-Q46} and \eqref{h47}, similar to \eqref{Q61} we get
     \begin{align} \label{h48}
   \psi_u''(s)\leq& \psi_u''(s)-6 \psi_u'(s) \nonumber\\
    =& -4 e^{2s}\int_{\mathbb{R}^3}|\nabla u|^2dx-\frac{5e^{s}}{4} \int_{\mathbb{R}^3} \int_{\mathbb{R}^3} \frac{|u(y)|^2|u(x)|^{2}} {|x-y|}dydx
+  (6-p\gamma_p)\gamma_p e^{p\gamma_ps}
\int_{\mathbb{R}^3}|u|^{p}dx   \nonumber\\
    \leq &  \left(-4+C_1(p)a^{\frac{8(p-3)}{3p-8}} \right)e^{2s} |\nabla u|_2^2\nonumber\\
    <& 0
   \end{align}
   when   $a<\left(\frac{4}{C_1(p)}\right)^{\frac{3p-8}{8(p-3)}}$.
 So, we obtain that $ \psi_u$   is    strictly     concave on  $(s_u, +\infty)$.
It remains to show that the map $u\in S_{a, r}\mapsto s_u\in \mathbb{R}$ is of class $C^1$.  Since $\pi(s,u):=\psi_u'(s)$ is of class $C^1 $, $\pi(s_u, u)=0 $ and $\partial_s\pi(s_u, u)=  \psi_u''(s_u)<0$, then the implicit function theorem gives that $u\in S_{a, r}\mapsto s_u\in \mathbb{R}$ is of class $C^1$.
Thus, we finish the proof.
 \end{proof}

\begin{lemma}\label{Q-Lem2.6}
 Let    $p\in(3, \frac{10}{3} )$ and $a\in(0, \overline{a})$.  For any $u\in \mathcal{P}_a$, there exists $\delta>0$ such that    $|\nabla u|_2^2> \delta$ and $m_a> 0$.
\end{lemma}
\begin{proof}
For any $u\in \mathcal{P}_a$, similar to \eqref{Q61} we have
\begin{align*}
0=&\int_{\mathbb{R}^3}|\nabla u|^2dx+\frac{1}{4}\int_{\mathbb{R}^3} \int_{\mathbb{R}^3} \frac{|u(y)|^2|u(x)|^{2}} {|x-y|}dydx
- \int_{\mathbb{R}^3}|u|^{6}dx
-\gamma_p\int_{\mathbb{R}^3}|u|^{p}dx\\
\geq& \left(1-C_2(p)a^{\frac{8(p-3)}{3p-8}} \right) |\nabla u|_2^2-S^{-3}|\nabla u|_2^6,
\end{align*}
where
 \begin{align}\label{k9}
    C_2(p):=\frac{3p-8}{2}\left[\gamma_p\right]^{\frac{2}{3p-8}}\left(\frac{2(10-3p)}{\widehat{C}}\right)^{\frac{10-3p}{3p-8}} \end{align}
which implies that there exists  $\delta>0$ satisfying  $|\nabla u|_2^2> \delta$  when $a<\big(\frac{1}{C_2(p)}\big)^{\frac{3p-8}{8(p-3)}}$. Moreover, using \eqref{j-Q46}, similar to \eqref{Q61} one has
\begin{align}
\mathcal{J}(u)&=\mathcal{J}(u) - \frac{1}{6} P(u)\nonumber\\
 &=\frac{1}{3} \int_{\mathbb{R}^3} |\nabla u|^2 dx + \frac{5}{24} \int_{\mathbb{R}^3} \int_{\mathbb{R}^3}  \frac{|u(y)|^2|u(x)|^2}{|x-y|} dy  dx +  \left( \frac{\gamma_p}{6} - \frac{1}{p} \right) \int_{\mathbb{R}^3} |u|^p dx\nonumber\\
 &\geq \frac{1}{3}|\nabla u|_2^2                                                                                                                 +\frac{5}{24} \widehat{C}  \frac{|u|_p^{\frac{2p}{10-3p}}}{ |\nabla u|_2^{\frac{2(3p-8)}{10-3p}} |u|_2^{\frac{8(p-3)}{10-3p}}}  - \frac{6-p\gamma_p}{6p} |u|_p^p\nonumber\\
 &\geq \left(\frac{1}{3} - C_3(p) a^{\frac{8(p-3)}{3p-8}}\right)|\nabla u|_2^2\label{k5}\\
  &\geq \left(\frac{1}{3} - C_3(p) a^{\frac{8(p-3)}{3p-8}}\right) \delta,\nonumber
\end{align}
where
\begin{align}\label{k11}
    C_3(p):=\frac{3p-8}{2}\left[\frac{6-p\gamma_p}{6p}\right]^{\frac{2}{3p-8}}\left(\frac{12(10-3p)}{5\widehat{C}}\right)^{\frac{10-3p}{3p-8}}. \end{align}
This infers that $m_a>\Big(\frac{1}{3} - C_3(p) a^{\frac{8(p-3)}{3p-8}}\Big) \delta>0$ when $a<\big(\frac{1}{3C_3(p)}\big)^{\frac{3p-8}{8(p-3)}}$.
\end{proof}

\begin{lemma}\label{Q-Lem2.4}
   Let   $p\in(3, \frac{10}{3} )$ and $a\in(0, \overline{a})$.   There is a constant $k > 0$ small enough such that
$$
0 < \sup_{\overline{A_k}} \mathcal{J} < m_{a} \quad \text{and} \quad u \in \overline{A_k} \implies \mathcal{J}(u), P(u) > 0,
$$
where $A_k = \{u \in S_{a, r} : |\nabla u|_2^2 < k\}$.
\end{lemma}
\begin{proof}
By \eqref{R5} and Lemma \ref{Q-Lem2.1}, one has
\begin{align}\label{R7}
  \mathcal{J}(u)
 &=\frac{1}{2}\int_{\mathbb{R}^3}|\nabla u|^2dx+\frac{1}{4} \int_{\mathbb{R}^3} \int_{\mathbb{R}^3} \frac{|u(y)|^2|u(x)|^{2}} {|x-y|}dydx
 -\frac{1}{6}\int_{\mathbb{R}^3}|u|^{6}dx-\frac{1}{p}\int_{\mathbb{R}^3}|u|^{p}dx\nonumber\\
 &\geq \frac{1}{2}                                                                                                                     |\nabla u|_2^2
+ \frac{1}{4} \overline{C}  \frac{|u|_p^{\frac{2p}{4-p}}}{ |\nabla u|_2^{\frac{2(2p-5)}{4-p}} |u|_2^{\frac{2(p-3)}{4-p}}}  -\frac{1}{p}|u|_p^p -\frac{1}{6}\mathcal{S}^{-3}|\nabla u|_2^6.
\end{align}
Let $K:= \frac{1}{4}\overline{C}  \Big(|\nabla u|_2^{\frac{2(2p-5)}{4-p}} |u|_2^{\frac{2(p-3)}{4-p}}\Big)^{-1} $.
 Consider the auxiliary function
 \begin{align*}
f(x)=\frac{1}{2}|\nabla u|_2^2 -\frac{1}{6}\mathcal{S}^{-3}|\nabla u|_2^6+ K x^{\frac{2}{4-p}}   -\frac{1}{p}x, \quad\quad x>0.
 \end{align*}
  By a simple calculation, one has
  \begin{align*}
f'(x) &= K \cdot \frac{2}{4-p} \cdot x^{\frac{p-2}{4-p}} - \frac{1}{p},\\
f''(x) &= K \cdot \frac{2}{4-p} \cdot \frac{p-2}{4-p} \cdot x^{\frac{p-2}{4-p}-1} > 0, \quad \quad \text{for all } x>0.
\end{align*}
Therefore, the function $f(x)$ has the unique global minimum at
\begin{align*}
\bar{x} = \left( \frac{4-p}{2pK} \right)^{\frac{4-p}{p-2}},
  \end{align*}
  and then
  \begin{align}\label{Q28}
  \min_{x>0}f(x)
  =f(\bar{x} )
  &= \frac{1}{2}                                                                                                                       |\nabla u|_2^2-\frac{1}{6}\mathcal{S}^{-3}|\nabla u|_2^6
   + K \cdot \left( \frac{4-p}{2pK} \right)^{\frac{2}{p-2}}
   - \frac{1}{p} \cdot \left( \frac{4-p}{2pK} \right)^{\frac{4-p}{p-2}}\nonumber\\
&= \frac{1}{2}                                                                                                                       |\nabla u|_2^2-\frac{1}{6}\mathcal{S}^{-3}|\nabla u|_2^6
   - \left( \frac{4-p}{2p} \right)^{\frac{4-p}{p-2}}
   \cdot \frac{p-2}{2p} \cdot K^{\frac{p-4}{p-2}} \nonumber\\
&= \frac{1}{2}|\nabla u|_2^2-\frac{1}{6}\mathcal{S}^{-3}|\nabla u|_2^6-C(p) a^{\frac{2(p-3)}{p-2}}|\nabla u|_2^{\frac{2(2p-5)}{p-2}},
    \end{align}
    where
    $$
    C(p)=\frac{p-2}{2p} \left( \frac{2(4-p)}{p\overline{C}} \right)^{\frac{4-p}{p-2}}.
    $$
  Back to \eqref{R5},   one gets
    \begin{align}\label{R6}
\mathcal{J} (u) &\geq  \frac{1}{2}                                                                                                                     |\nabla u|_2^2
+ \frac{1}{4} \overline{C} \frac{|u|_p^{\frac{2p}{4-p}}}{ |\nabla u|_2^{\frac{2(2p-5)}{4-p}} |u|_2^{\frac{2(p-3)}{4-p}}}  -\frac{1}{p}|u|_p^p -\frac{1}{6}\mathcal{S}^{-3}|\nabla u|_2^6\nonumber\\
 &\geq \frac{1}{2}|\nabla u|_2^2-\frac{1}{6}\mathcal{S}^{-3}|\nabla u|_2^6-C(p) a^{\frac{2(p-3)}{p-2}}|\nabla u|_2^{\frac{2(2p-5)}{p-2}}\\
 &\geq  |\nabla u|_2^2 \left(\frac{1}{2}-C(p) a^{\frac{2(p-3)}{p-2}} |\nabla u|_2^{\frac{2(p-3)}{p-2}}-\frac{1}{6}\mathcal{S}^{-3}|\nabla u|_2^4\right).\nonumber
\end{align}
Since $\frac{2(p-3)}{p-2}>0$,  there exists $k > 0$ small enough such that
$\mathcal{J}(u)>0$ for any $u\in\overline{A_k}$.  Moreover, similar to \eqref{R6} we have
\begin{align*}
P(u)=&\int_{\mathbb{R}^3}|\nabla u|^2dx+\frac{1}{4}\int_{\mathbb{R}^3} \int_{\mathbb{R}^3} \frac{|u(y)|^2|u(x)|^{2}} {|x-y|}dydx
- \int_{\mathbb{R}^3}|u|^{6}dx
-\gamma_p\int_{\mathbb{R}^3}|u|^{p}dx\\
\geq& \left(1-C(p)\cdot (p\gamma_p)^{\frac{2}{p-2}} a^{\frac{2(p-3)}{p-2}} |\nabla u|_2^{\frac{2(p-3)}{p-2}}- \mathcal{S}^{-3}|\nabla u|_2^4\right)|\nabla u|_2^2,
\end{align*}
which show that  $P(u)>0$ for any $u\in\overline{A_k}$ with $k > 0$ small enough.
Besides,  choosing $k$ sufficient small, from Lemma \ref{Q-Lem2.6} we
have, for all $u\in\overline{A_k}$
\begin{align*}
 \mathcal{J}(u)
 &=\frac{1}{2}\int_{\mathbb{R}^3}|\nabla u|^2dx+\frac{1}{4} \int_{\mathbb{R}^3} \int_{\mathbb{R}^3} \frac{|u(y)|^2|u(x)|^{2}} {|x-y|}dydx
 -\frac{1}{6}\int_{\mathbb{R}^3}|u|^{6}dx-\frac{1}{p}\int_{\mathbb{R}^3}|u|^{p}dx\\
 &\leq \frac{1}{2} |\nabla u|_2^2+Ca^3 |\nabla u|_2\\
 &<m_a.
\end{align*}
Hence, we complete the proof of this lemma.
\end{proof}

\begin{lemma}\label{Lem4.1}
  Let    $p\in(3, \frac{10}{3} )$. If $ u \in \mathcal{P}_{a}$ is a critical point for $\mathcal{J}|_{\mathcal{P}_{a} }$ and $\left(\mathcal{P}_{a} \right)^{ 0}=\emptyset$, then    $u$ is a critical point for $\mathcal{J}|_{S_{a, r} } $.
\end{lemma}

\begin{proof}
We claim that $dP(u)$ and $dG(u)$  are linearly independent, where $G(u)=|u|_2^2-a^2$. Indeed, let
\begin{align}\label{h5}
\nu_1P'(u)+\nu_2u=0,
\end{align}
where $\nu_1, \nu_2\in \mathbb{R}$. Then one has
$$
-2\nu_1(\Delta u)+\nu_2 u+\nu_1\left(|x|^{-1}\ast|u|^2\right)u-6\nu_1|u|^{4}u-\nu_1p\gamma_p|u|^{p-2}u=0.
$$
 The corresponding Nehari identity and Poho\v{z}aev identity are given by
\begin{align}\label{h21}
  \nu_1\Big[2|\nabla u|_2^2+  \int_{\mathbb{R}^3}  \int_{\mathbb{R}^3} \frac{|u(y)|^2|u(x)|^2}{|x-y|} dydx -6 |u|_{6}^{6}-p\gamma_p|u|_{p}^{p}\Big]+\nu_2|u|_2^2=0
\end{align}
and
\begin{align}\label{h22}
 \nu_1 \Big[ |\nabla u|_2^2+ \frac{5}{4} \int_{\mathbb{R}^3}  \int_{\mathbb{R}^3} \frac{|u(y)|^2|u(x)|^2}{|x-y|} dydx     -3 |u|_{6}^{6}-3 \gamma_p|u|_{p}^{p}\Big]+\frac{3}{2}\nu_2|u|_2^2=0.
\end{align}
Combining \eqref{h21} with \eqref{h22}, one gets
\begin{align}\label{h7}
\nu_1 \Big(2|\nabla u|_2^2+ \frac{1}{4}  \int_{\mathbb{R}^3}  \int_{\mathbb{R}^3} \frac{|u(y)|^2|u(x)|^2}{|x-y|} dydx - 6|u|_{6}^{6}- p \gamma_p^2 |u|_{p}^{p}\Big )=0.
\end{align}
Since $\left(\mathcal{P}_{a} \right)^{ 0}=\emptyset$, then
\begin{align}\label{h6}
2 \int_{\mathbb{R}^3} |\nabla u|^2 dx + \frac{1}{4}\int_{\mathbb{R}^3} \int_{\mathbb{R}^3} \frac{|u(y)|^2|u(x)|^{2}} {|x-y|}dydx - 6 |u|_6^6 -  p\gamma_p^2 |u|_p^p\neq0,
\end{align}
 which and \eqref{h7} show that $\nu_1=0$. From \eqref{h5}, one infers $\nu_2=0$. Thus, the claim holds. Besides,
if $ u \in \mathcal{P}_{a}$ is a critical point for $\mathcal{J}|_{\mathcal{P}_{a} }$,   by the Lagrange multipliers rule,  there exist $ \lambda, \nu\in \mathbb{R}$ such that
$$
\mathcal{J}'(u)+\lambda u+\nu P'(u)=0,
$$
namely,
$$
(1+2\nu)(-\Delta u)+(1+\nu ) \left(|x|^{-1}\ast|u|^2\right)u +\lambda u-(1+6 \nu )|u|^{4}u-(1+ \nu p\gamma_p  )|u|^{p-2}u=0 \ \ \ \text{in} \ \mathbb{R}^3.
$$
The corresponding Poho\v{z}aev identity is given by
\begin{align*}
(1+2\nu) |\nabla u|_2^2+\frac{1}{4}(1+\nu )\int_{\mathbb{R}^3}  \int_{\mathbb{R}^3} \frac{|u(y)|^2|u(x)|^2}{|x-y|} dydx  -(1+6 \nu ) |u|_{6}^{6}-\gamma_p(1+ \nu p\gamma_p)|u|_{p}^{p}=0,
\end{align*}
which and  $P(u)=0$ imply that
$$
\nu \big(2|\nabla u|_2^2+ \frac{1}{4}\int_{\mathbb{R}^3} \int_{\mathbb{R}^3} \frac{|u(y)|^2|u(x)|^{2}} {|x-y|}dydx - 6|u|_{6}^{6}- p \gamma_p^2 |u|_{p}^{p}=0.
$$
Then,  we get from \eqref{h6} that $\nu=0$. Hence, we complete the proof.
\end{proof}

\section{Proof of Theorem \ref{K-TH1} }\label{sec3}
To prepare for the existence of standing wave solutions to Eq. \eqref{L1}, we first construct a Palais-Smale sequence that asymptotically satisfies the Poho\v{z}aev identity.
\begin{lemma}\label{R-lem1}
Let $3 <p < \frac{10}{3}$ and $a\in(0, \overline{a})$.   Then there exists a sequence $\{u_n\} \subset \mathcal{S}_{a, r} $ such that, as $n\rightarrow\infty$,
$$
\mathcal{J}(u_n) \to m_a > 0, \quad \mathcal{J}|'_{\mathcal{S}_{a,r}}(u_n) \to 0 \quad \text{and} \quad P(u_n) \to 0.
$$
\end{lemma}
\begin{proof}
 Similar to \cite{1997jean},
considering the augmented functional  $\widetilde{\mathcal{J}}(s, u): \mathbb{R}\times \mathcal{H}_r\rightarrow \mathbb{R}$ defined by
\begin{align}\label{Q11}
 \widetilde{\mathcal{J}}(s, u ):=\mathcal{J}(s\star u)
 &=\frac{e^{2s}}{2}\int_{\mathbb{R}^3}|\nabla u|^2dx+\frac{ e^{s}}{4} \int_{\mathbb{R}^3} \int_{\mathbb{R}^3} \frac{|u(y)|^2|u(x)|^{2}} {|x-y|}dydx \nonumber\\
 &\quad-\frac{e^{6s}}{6}\int_{\mathbb{R}^3}|u|^{6}dx
 -\frac{e^{p\gamma_ps}}{p}
 \int_{\mathbb{R}^3}|u|^{p}dx,
\end{align}
 notice that $\widetilde{\mathcal{J}}$ is of class $C^1$. Let $k>0$ be defined by Lemma \ref{Q-Lem2.4}.
 Denoting by $\mathcal{J}^c$ the closed sub-level set $\{u \in S_{a, r}: \mathcal{J}(u)\leq c\}$, we introduce the minimax class
  \begin{align*}
\Gamma:=\left\{\gamma=(\alpha,\beta)\in C([0,1], \mathbb{R}\times S_{a, r}): \gamma(0)\in (0, \overline{A}_k), \gamma(1)\in (0,  \mathcal{J}^0) \right\},
  \end{align*}
with associated minimax level
  \begin{align*}
 c_{a}:=\inf_{\gamma\in \Gamma}\max_{(s,u)\in \gamma([0,1])} \widetilde{\mathcal{J}}(s,u).
  \end{align*}
 Let $u\in S_{a, r}$. Since $|\nabla (s\star u)|_2\rightarrow 0^+$ as $s\rightarrow -\infty$ and $  \mathcal{J}(s\star u)\rightarrow -\infty$ as  $s\rightarrow +\infty$, then there exist $s_0\ll-1$ and $s_1\gg1$ such that
  \begin{align}\label{h28}
  \gamma_u: \tau\in[0,1]\mapsto (0, ((1-\tau)s_0+\tau s_1)\star u)\in \mathbb{R}\times S_{a, r}
   \end{align}
  is a path in $\Gamma$. Then $c_{a}$ is a real number. Now, for any $\gamma=(\alpha,\beta)\in \Gamma$, let us consider the function
  $$
  P\circ\gamma:\tau\in[0,1] \mapsto P(\alpha(\tau)\star \beta(\tau))\in \mathbb{R}.
  $$
By Lemma \ref{Q-Lem2.4}, we have  $P\circ\gamma(0) =P(\beta(0)) >0$ since $\beta(0)\in\overline{ A}_k$.   Now we claim that
$$
P\circ\gamma(1) =P(\beta(1)) <0.
$$
In fact, it follows from  \eqref{k8} and \eqref{h20} that
 $ \psi_u(s)>0$ for $s\rightarrow-\infty$ and $ \psi_u'(s)>0$ for $s\rightarrow-\infty$.
We know from Lemma \ref{Q-Lem2.3} that $ \psi_{\beta(1)}(s)>0$ for $s\leq s_{\beta(1)}$.
 Moreover, since $\beta(1)\in \mathcal{J}^0$,  we have
 $\psi_{\beta(1)}(0)=\mathcal{J}(\beta(1))\leq 0,$
 which implies that $s_{\beta(1)}<0$. Thus, from Lemma \ref{Q-Lem2.3} we get the claim. Furthermore, the map $\tau\mapsto  \alpha(\tau)\star\beta(\tau)$ is continuous from $[0, 1]$ to $\mathcal{H}_r$, and hence we deduce that there exists  $\tau_\gamma\in(0, 1)$ such that $P\circ\gamma(\tau_\gamma ) =0$, namely $\alpha(\tau_\gamma)\star\beta(\tau_\gamma)\in \mathcal{P}_{a}$; this implies that
$$
\max_{(s,u)\in\gamma([0,1])} \widetilde{\mathcal{J}}(s,u)\geq  \widetilde{\mathcal{J}}( \gamma(\tau_\gamma))=\widetilde{\mathcal{J}}(\alpha(\tau_\gamma), \beta(\tau_\gamma) )= \mathcal{J}(\alpha(\tau_\gamma)\star\beta(\tau_\gamma) )\geq \inf_{u\in\mathcal{P}_{a} } \mathcal{J}=m_{a}.
$$
Hence, $ c_{a}\geq m_{a}$. In addition, if $u\in \mathcal{P}_{a}  $, then $\gamma_u$  defined in \eqref{h28} is a path in $\Gamma$ with
$$
\mathcal{J}(u)\geq \widetilde{\mathcal{J}}(\gamma_u)\geq c_{a},
$$
whence the reverse inequality $m_{a} \geq c_{a}$ follows.  Combining this with Lemma \ref{Q-Lem2.4}, we infer that
$$
c_{a}=m_{a}>\sup_{( \overline{A}_k\cup \mathcal{J}^0)\cap S_{a, r}}\mathcal{J}=\sup_{( (0,\overline{A}_k)\cup (0,\mathcal{J}^0))\cap (\mathbb{R}\times S_{a, r})} \widetilde{\mathcal{J}}.
$$
Using the terminology in \cite[Section 5]{1993book}, this means that $\{\gamma([0, 1]) :\gamma\in \Gamma\}$ is a homotopy stable family of compact subsets of $\mathbb{R} \times S_{a, r}$ with extended closed boundary $(0, \overline{A}_k)\cup (0,  \mathcal{J}^0)$, and that the superlevel set $\{ \mathcal{J}\geq c_{a}\}$ is a dual set for $\Gamma$.  Therefore, applying \cite[Theorem 3.2]{1993book},  choosing any minimizing sequence ${\gamma_n=(\alpha_n, \beta_n)}\subset \Gamma$ for $\widetilde{\mathcal{J}}|_{\mathbb{R}\times S_{a, r}}$  at the level   $c_{a} $ with the property that $\alpha_n=0$ and $\beta_n(\tau) \geq0$ a.e. in $\mathbb{R}^3$ for every $\tau\in [0, 1]$, there exists a   sequence $\{(s_n, w_n)\}\subset \mathbb{R} \times S_{a, r}$  such that, as  $n\rightarrow\infty$,
\begin{enumerate}
    \item[$(i)$] $\widetilde{\mathcal{J}}(s_n, w_n) \to c_a$;
    \item[$(ii)$] $\widetilde{\mathcal{J}}'|_{\mathbb{R} \times S_{a, r}}(s_n, w_n) \to 0$;
    \item[$(iii)$] $\operatorname{dist}((s_n, w_n), (0, \beta_n(\tau))) \to 0$.
\end{enumerate}
Let $u_n := s_n \star w_n = e^{\frac{3 s_n}{2}} w_n(e^{s_n} x)$. one obtains from $(i)$ that
\begin{align}\label{Q36}
\lim_{n \to \infty} \mathcal{J}(u_n) = \lim_{n \to \infty} \mathcal{J}(s_n \star  w_n) = \lim_{n \to \infty} \widetilde{\mathcal{J}}(s_n, w_n) = c_a.
\end{align}
In virtue of $(ii)$ and \eqref{Q11}, we deduce
$$
\partial_s \widetilde{ \mathcal{J}}(s_n, w_n) = \big\langle \widetilde{\mathcal{J}}'(s_n, w_n), (1, 0) \big\rangle \to 0 \quad\quad \text{as } n \to \infty,
$$
where
\begin{align*}
\partial_s \widetilde{ \mathcal{J}}(s_n, w_n)
= &e^{2s_n} \int_{\mathbb{R}^3} |\nabla w_n|^2 dx +\frac{e^{s_n}}{4}\int_{\mathbb{R}^3} \int_{\mathbb{R}^3} \frac{|w_n(y)|^2|w_n(x)|^2}{|x-y|}dy dx \\
&- e^{6 s_n} \int_{\mathbb{R}^3} |w_n|^{6} dx -  \gamma_p e^{3(\frac{p}{2} - 1)s_n} \int_{\mathbb{R}^3} |w_n|^p dx.
\end{align*}
Thus, from \eqref{Q30} we get
\begin{align}\label{Q35}
P(u_n) = P(s_n\star w_n) = \partial_s \widetilde{\mathcal{J}}(s_n, w_n) \to 0 \quad \quad \text{as}\ n \to \infty.
\end{align}
Besides, setting $h_n \in T_{u_n} := \left\{ z \in \mathcal{H}_r: \langle z, u_n \rangle_{L^2} = 0 \right\}$, by simple calculation, one has
\begin{align*}
\langle \mathcal{J}'(u_n), h_n \rangle
&= \int_{\mathbb{R}^3} \nabla u_n \cdot \nabla h_ndx  + \int_{\mathbb{R}^3}\int_{\mathbb{R}^3} \frac{|u_n(y)|^2}{|x-y|}u_n(x) h_n(x)dydx\\
  &\quad - \int_{\mathbb{R}^3} |u_n|^{4} u_n h_n dx
   -  \int_{\mathbb{R}^3} |u_n|^{p-2} u_n h_n dx \\
&= e^{2s_n} \int_{\mathbb{R}^3} e^{-\frac{5}{2} s_n} \nabla w_n \cdot \nabla h_n(e^{-s_n}x) dx
   + e^{s_n}\int_{\mathbb{R}^3} \int_{\mathbb{R}^3} \frac{|w_n(y)|^2}{|x-y|}w_n(x)h_n(e^{-s_n }x)e^{-\frac{3}{2}s_n} dydx \\
&\quad - e^{6 s_n} \int_{\mathbb{R}^3} |w_n|^{4} w_n e^{-\frac{3 s_n}{2}} h_n(e^{-s_n}x) dx
   -  e^{3 s_n(\frac{p}{2}-1) } \int_{\mathbb{R}^3} |w_n|^{p-2} w_n e^{-\frac{3s_n}{2}} h_n(e^{-s_n}x) dx \\
&= e^{2s_n} \int_{\mathbb{R}^3} \nabla w_n \cdot \nabla \widetilde{h}_n dx
   + e^{s_n} \int_{\mathbb{R}^3}  \int_{\mathbb{R}^3} \frac{|w_n(y)|^2}{|x-y|} w_n(x) \widetilde{h}_n(x) dydx \\
&\quad - e^{6 s_n} \int_{\mathbb{R}^3} |w_n|^{4} w_n \widetilde{h}_n dx
   -  e^{3(\frac{p}{2}-1)  s_n} \int_{\mathbb{R}^3} |w_n|^{p-2} w_n \widetilde{h}_n dx,
\end{align*}
which infers that
\begin{align}\label{Q32}
\langle \mathcal{J}'(u_n), h_n \rangle = \big\langle \widetilde{\mathcal{J}}'(s_n, w_n), (0, \widetilde{h}_n) \big\rangle,
\end{align}
where $\widetilde{h}_n(x) = e^{-\frac{3s_n}{2}} h_n(e^{-s_n}x)$. Now, we claim that
\begin{align}\label{Q33}
(0, \widetilde{h}_n) \in \widetilde{T}_{(s_n, w_n)} := \left\{ (z_1, z_2) \in \mathbb{R} \times \mathcal{H}_r: \langle w_n, z_2 \rangle_{L^2} = 0 \right\}.
\end{align}
Indeed,
\begin{align*}
(0, \widetilde{h}_n) \in \widetilde{T}_{(s_n, w_n)}
& \Leftrightarrow \langle w_n, \widetilde{h}_n \rangle_{L^2} = 0 \\
& \Leftrightarrow \int_{\mathbb{R}^3} w_n e^{-\frac{3s_n}{2}} h_n(e^{-s_n}x) dx = 0 \\
& \Leftrightarrow \int_{\mathbb{R}^3} e^{\frac{3 s_n}{2}} w_n(e^{s_n}x) h_n(x) dx = 0 \\
& \Leftrightarrow \langle u_n, h_n \rangle_{L^2} = 0 \\
& \Leftrightarrow h_n \in T_{u_n}.
\end{align*}
Furthermore, using the last item $(iii)$, we conclude that $s_n\rightarrow 0$ as $n\rightarrow\infty$. 
Thus, in view of $(ii)$, \eqref{Q32}, \eqref{Q33}, and the fact that $\|\widetilde{h}_n\| \leq C \|h_n\|$ one has, as $n \to \infty$,
$$
\langle \mathcal{J}'(u_n), h_n \rangle \to 0 \quad\quad \text{for all } h_n \in T_{u_n},
$$
which yields  that
\begin{align}\label{Q34}
\mathcal{J}'|_{S_{a, r}}(u_n) \to 0 \quad \quad \text{as} \ n \to \infty.
\end{align}
Thus, from \eqref{Q36}, \eqref{Q35} and \eqref{Q34}, we find that $\{u_n\} \subset S_{a, r}$ is a Palais-Smale sequence for $\mathcal{J} $ at the level $c_a = m_a$ satisfying
$$
\mathcal{J}(u_n) \to m_a, \quad \mathcal{J}'|_{S_{a, r}}(u_n) \to 0 \quad  \text{and}  \quad P(u_n) \to 0 \quad \text{as } n \to \infty.
$$
Hence, we complete the proof.
\end{proof}

Now, we derive an upper bound estimate for $m_a$, which plays an essential role in proving compactness.

\begin{lemma}\label{R-lem2}
Let $3<p<\frac{10}{3}$. Then there exist constants $a_*>0$    and
$\eta\in\big(4,\frac{2(5p-12)}{3p-8}\big)$ such that, for any $a\in(0, a_*)$,
\begin{align*}
m_a\leq \frac{1}{3}\mathcal{S}^{\frac{3}{2}}-c_pa^{\eta},
\end{align*}
where $c_p>0$ is a constant, independent of $a$.
\end{lemma}

\begin{proof}
Denote $W_{\varepsilon}(\cdot)=\chi(\cdot)U_{\varepsilon}(\cdot)$, where
\begin{align*}
U_{\varepsilon}(x)=3^{\frac14}\left(\frac{\varepsilon}{\varepsilon^2+|x|^2}\right)^{\frac12}.
\end{align*}
Take a nonnegative radial cut-off function $\chi\in C_0^\infty(\mathbb R^3)$ such that
$\chi(x)=1$ for $|x|\leq1$ and $\chi(x)=0$ for $|x|\geq2$.
By the well-known computations in \cite{1983-CPAM-Brezis}, we have, as $\varepsilon\to0^+$,
\begin{align}\label{I8}
|\nabla W_{\varepsilon}|_2^2=\mathcal{S}^{\frac32}+O(\varepsilon),\qquad
|W_{\varepsilon}|_6^6=\mathcal{S}^{\frac32}+O(\varepsilon^3)
\end{align}
and
\begin{align}\label{I1}| W_{\varepsilon}|_{p}^{p}\sim
\begin{cases}
\displaystyle \varepsilon^{\frac{p}{2}},  &\text{if}\ p\in[2,3);\\
\displaystyle \varepsilon^{\frac{3}{2}}\ln \frac{1}{\varepsilon}, &\text{if}\ p=3;\\
\displaystyle \varepsilon^{3-\frac{p}{2}}, &\text{if}\ p\in(3, 6).
\end{cases}
\end{align}
Moreover, by Lemma \ref{gle4}, a direct calculation yields
\begin{align*}
\Lambda(W_{\varepsilon}):=\int_{\mathbb R^3}\int_{\mathbb R^3}
\frac{|W_{\varepsilon}(x)|^2|W_{\varepsilon}(y)|^2}{|x-y|}dxdy \leq C|W_{\varepsilon}|_{\frac{12}{5}}^4
\leq C\varepsilon^2.
\end{align*}
Let $v_{\varepsilon}=aW_{\varepsilon}/|W_{\varepsilon}|_2$. Then
$v_{\varepsilon}\in S_{a,r}$. By Lemma \ref{Q-Lem2.3}, there is a unique
$t_{\varepsilon}\in\mathbb R$ such that
$t_{\varepsilon}\star v_{\varepsilon}\in\mathcal P_a$, and
\begin{align}\label{Q55}
m_a\leq\mathcal J(t_{\varepsilon}\star v_{\varepsilon})
=\max_{t\in\mathbb R}\mathcal J(t\star v_{\varepsilon}).
\end{align}
Setting
\begin{align*}
r=\frac{ae^t}{|W_{\varepsilon}|_2}>0,
\end{align*}
 we obtain
\begin{align}\label{Q53}
\mathcal J(t\star v_{\varepsilon})
&=\frac{A_{\varepsilon}}2r^2-\frac{B_{\varepsilon}}6r^6
+\frac{d_{a,\varepsilon}}4r
-\frac{b_{a,\varepsilon}}p r^{\frac{3(p-2)}2}
=:F_{a,\varepsilon}(r),
\end{align}
where
\begin{align*}
A_{\varepsilon}&=|\nabla W_{\varepsilon}|_2^2,\qquad
B_{\varepsilon}=|W_{\varepsilon}|_6^6,\\
d_{a,\varepsilon}&=\frac{a^3\Lambda(W_{\varepsilon})}{|W_{\varepsilon}|_2^3},\qquad
b_{a,\varepsilon}=\frac{a^{3-\frac p2}|W_{\varepsilon}|_p^p}
{|W_{\varepsilon}|_2^{3-\frac p2}}.
\end{align*}
In view of the above estimates, there exist constants $c,C>0$, independent of
$a$ and $\varepsilon$, such that
\begin{align*}
&A_{\varepsilon}=\mathcal S^{\frac32}+O(\varepsilon), \quad\quad\quad
B_{\varepsilon}=\mathcal S^{\frac32}+O(\varepsilon^3),\\
&0\leq d_{a,\varepsilon}\leq Ca^3\varepsilon^{\frac12}, \quad\quad
ca^{3-\frac p2}\varepsilon^{\frac{6-p}{4}}
\leq b_{a,\varepsilon}
\leq Ca^{3-\frac p2}\varepsilon^{\frac{6-p}{4}}.
\end{align*}
Next, we choose $\varepsilon$ depending on $a$. Since
$3<p<\frac{10}{3}$, a direct calculation shows that
\begin{align*}
4<\frac{2(6-p)}{p-2}
<\min\left\{\frac{2p}{4-p},\frac{2(5p-12)}{3p-8}\right\}.
\end{align*}
We may therefore choose $\kappa$ sufficiently close to
$\frac{2(6-p)}{p-2}$ from above such that
\begin{align*}
\frac{2(6-p)}{p-2}<\kappa<\frac{2p}{4-p},
\end{align*}
and put
\begin{align*}
\eta=3-\frac p2+\frac{6-p}{4}\kappa.
\end{align*}
It follows that
\begin{align*}
4<\eta<\frac{2(5p-12)}{3p-8},\qquad
\eta<\kappa,\qquad
\eta<3+\frac{\kappa}{2}.
\end{align*}
From now on, let $\varepsilon=a^\kappa$ and write
$F_a(r):=F_{a,a^\kappa}(r)$. For $a\to0^+$ we have
\begin{align*}
&A_{a^\kappa}=\mathcal S^{\frac32}+O(a^\kappa),\qquad
B_{a^\kappa}=\mathcal S^{\frac32}+O(a^{3\kappa}),\\
&0\leq d_{a,a^\kappa}\leq Ca^{3+\frac{\kappa}{2}},\qquad
ca^\eta\leq b_{a,a^\kappa}\leq Ca^\eta.
\end{align*}
We now estimate the global maximum of $F_a$. Denote $T=\mathcal S^{\frac32}$ and
\begin{align*}
F_0(r)=T\left(\frac{r^2}{2}-\frac{r^6}{6}\right),\qquad r\geq0.
\end{align*}
Clearly, $F_0'(r)=Tr(1-r^4)$, and hence $r=1$ is the unique global
maximum point of $F_0$, with $F_0(1)=T/3$. The preceding estimates imply
that $F_a\to F_0$ uniformly on every bounded interval of $[0,+\infty)$.
For all sufficiently small $a>0$, we have
\begin{align*}
F_a(r)\leq Tr^2+\frac14r-\frac{T}{12}r^6,\qquad r\geq0.
\end{align*}
Thus, there exists $R>2$, independent of $a$, such that $F_a(r)<0$ for
all $r\geq R$, whereas $F_a(1)\to T/3>0$. In particular, every $F_a$
attains its global maximum at some $r_a\in(0,R)$.
We claim that
\begin{align*}
r_a\to1\qquad\text{as }a\to0^+.
\end{align*}
Indeed, if the claim does not hold, one can choose $a_n\to0^+$ and
corresponding maximum points $r_{a_n}\to r_0\in[0,R]$ with $r_0\neq1$.
By the uniform convergence of $F_{a_n}$ to $F_0$ on $[0,R]$,
\begin{align*}
F_0(r_0)=\lim_{n\to\infty}F_{a_n}(r_{a_n})
\geq\lim_{n\to\infty}F_{a_n}(1)=F_0(1),
\end{align*}
which contradicts the uniqueness of the maximum point of $F_0$.
Thus, for $a>0$ sufficiently small,
\begin{align*}
\frac12\leq r_a\leq\frac32.
\end{align*}
Combining this fact with the estimates for $A_{a^\kappa}$,
$B_{a^\kappa}$, $d_{a,a^\kappa}$ and $b_{a,a^\kappa}$, we deduce that
\begin{align*}
\max_{r>0}F_a(r)
&=F_a(r_a)\\
&\leq F_0(r_a)+Ca^\kappa+Ca^{3\kappa}
+Ca^{3+\frac\kappa2}-ca^\eta\\
&\leq\frac13\mathcal S^{\frac32}+Ca^\kappa+Ca^{3\kappa}
+Ca^{3+\frac\kappa2}-ca^\eta.
\end{align*}
Since $\eta<\kappa$ and $\eta<3+\frac\kappa2$,  there exist $a_*>0$ and  $c_p>0$ such that
\begin{align*}
\max_{r>0}F_a(r)\leq\frac13\mathcal S^{\frac32}-c_pa^\eta
\qquad\text{for all }a\in(0, a_*).
\end{align*}
Finally, using \eqref{Q55} and \eqref{Q53}, one has
\begin{align*}
m_a\leq\max_{t\in\mathbb R}\mathcal J(t\star v_{a^\kappa})
=\max_{r>0}F_a(r)
\leq\frac13\mathcal S^{\frac32}-c_pa^\eta \qquad\text{for all }a\in(0, a_*).
\end{align*}
We complete the proof.
\end{proof}

 \begin{lemma}\label{R-lem6}
  Assume that $3 <p < \frac{10}{3}$ and $a\in(0, \overline{a_*})$, where $\overline{a_*}>0$ is sufficiently small.  Let $\{u_n\} \subset S_{a,r}$ be a Palais-Smale sequence at level $m_a \in (0, \frac{1}{3} \mathcal{S}^{\frac{3}{2}}-c_pa^\eta]$ such that $P(u_n) \to 0$ as $n \to \infty$. Then there exists $u \in S_{a}$ such that $u_n \rightarrow u$ in $\mathcal{H}_r$ and $(\lambda, u) \in \mathbb{R}^+ \times \mathcal{H}_r$ solves Eq. \eqref{L1}.
\end{lemma}

\begin{proof}
  We first  prove that $\{u_n\}$ is bounded in $\mathcal{H}_r$. Since $P(u_n) \to 0$, by the Gagliardo-Nirenberg inequality \eqref{k10}, one has
\begin{align*}
m_a + o_n(1) &= \mathcal{J}(u_n) - \frac{1}{6} P(u_n) \\
&= \frac{1}{3} \int_{\mathbb{R}^3} |\nabla u_n|^2 dx + \frac{5 }{24} \int_{\mathbb{R}^3}\int_{\mathbb{R}^3} \frac{|u_n(y)|^2}{|x-y|} dy |u_n|^2 dx +   \left( \frac{\gamma_p}{6} - \frac{1}{p} \right) \int_{\mathbb{R}^3} |u_n|^p dx \\
&\ge \frac{1}{3} \int_{\mathbb{R}^3} |\nabla u_n|^2 dx -   \frac{6 - p\gamma_p}{6p} C_p^p a^{(1-\gamma_p)p} |\nabla u_n|_2^{p\gamma_p},
\end{align*}
which shows that $\{u_n\}$ is bounded in  $\mathcal{H}_r$ because of $p\gamma_p < 2$.
Hence, there exists $u \in \mathcal{H}_r$ such that, up to subsequences,
$$
u_n \rightharpoonup u \ \text{in} \ \mathcal{H}_r; \quad u_n \to u \ \ \text{in} \ L^t(\mathbb{R}^3) \ \text{with} \ t \in (2,6); \quad u_n \to u \text{ a.e in } \mathbb{R}^3.
$$
Recall that $\mathcal{J}'|_{S_{a, r}}(u_n) \to 0$, then there exists $\lambda_n \in \mathbb{R}$ satisfying
\begin{align}\label{Q29}
-\Delta u_n + \lambda_n u_n +   \int_{\mathbb{R}^3} \frac{|u_n(y)|^2}{|x-y|} dy u_n
= |u_n|^4 u_n +   |u_n|^{p-2} u_n + o_n(1) \quad\quad \text{in } \mathcal{H}_r^{-1}.
\end{align}
Multiplying the above equation by $u_n$ and integrating in $\mathbb{R}^3$,
\begin{align}\label{Q30}
\lambda_n a^2 = -|\nabla u_n|_2^2 -   \int_{\mathbb{R}^3}  \int_{\mathbb{R}^3} \frac{|u_n(y)|^2}{|x-y|} dy |u_n|^2 dx
+ \int_{\mathbb{R}^3} |u_n|^6 dx +   \int_{\mathbb{R}^3} |u_n|^p dx + o_n(1).
\end{align}
Since $\{u_n\} \subset \mathcal{H}_r$ is bounded, then it follows from \eqref{Q30} that $\{\lambda_n\}$ is bounded in $\mathbb{R}$. Hence, there is $\lambda \in \mathbb{R}$ satisfying $\lambda_n \to \lambda$. Then we have
\begin{align}\label{Q37}
-\Delta u + \lambda u +   \int_{\mathbb{R}^3} \frac{|u(y)|^2}{|x-y|} dy u
= |u|^4 u +   |u|^{p-2} u \quad\quad \text{in } \mathcal{H}_r^{-1}.
\end{align}
Clearly, $u \neq 0$. If not, using $P(u_n) \to 0$ as $n\rightarrow\infty$, Lemma \ref{gle4} and the fact that $\mathcal{H}_r \hookrightarrow L^t(\mathbb{R}^3)$ is compact for $t \in (2,6)$, one gets
\begin{align}\label{Q38}
\mathcal{S} |u_n|_6^2 \leq |\nabla u_n|_2^2 \leq |u_n|_6^6 + o_n(1),
\end{align}
which produces two cases: either (i) $\limsup_{n \to \infty} |u_n|_6^6 = 0$ or (ii) $\limsup_{n \to \infty} |u_n|_6^6 = l > 0$. If (i) holds, by \eqref{Q38} we have $|\nabla u_n|_2^2 \to 0$, which implies $\mathcal{J}(u_n) \to 0$ contradicting $\mathcal{J}(u_n) \to m_a > 0$. If (ii) holds, from \eqref{Q38} one has $|u_n|_6^6 \geq \mathcal{S}^{\frac{3}{2}} + o_n(1)$. Since $\lim_{n \to \infty} |u_n|_p^p = 0$ and $\lim_{n \to \infty} P(u_n) = 0$, then
$$
|\nabla u_n|_2^2 = |u_n|_6^6 + o_n(1) \geq \mathcal{S}^{\frac{3}{2}} + o_n(1).
$$
Hence,
$$
\mathcal{J}(u_n) = \frac{1}{3} |\nabla u_n|_2^2 + o_n(1) \geq \frac{1}{3} \mathcal{S}^{\frac{3}{2}} + o_n(1),
$$
which implies $m_a \geq \frac{1}{3} \mathcal{S}^{\frac{3}{2}}$, a contradiction due to $m_a <  \frac{1}{3} \mathcal{S}^{\frac{3}{2}}$. Thus, $u \neq 0$. Besides, we claim that
\begin{align}\label{h85}
\lambda > 0.
\end{align}
Since $(\lambda, u) \in \mathbb{R} \times \mathcal{H}_r$ satisfies \eqref{Q37}, then the following equalities hold:
\begin{align}
|\nabla u|_2^2 + \lambda \int_{\mathbb{R}^3} |u|^2 dx
+  \int_{\mathbb{R}^3}\int_{\mathbb{R}^3} \frac{|u(y)|^2}{|x-y|} dy |u|^2 dx
- \int_{\mathbb{R}^3} |u|^6 dx -   \int_{\mathbb{R}^3} |u|^p dx = 0,\label{Q40} \\
|\nabla u|_2^2 + \frac{1}{4} \int_{\mathbb{R}^3}\int_{\mathbb{R}^3} \frac{|u(y)|^2}{|x-y|} dy |u|^2 dx
- \int_{\mathbb{R}^3} |u|^6 dx -   \gamma_p \int_{\mathbb{R}^3} |u|^p dx = 0.
\label{Q41}
\end{align}
Assume on the contrary that $\lambda \leq 0$. By combining \eqref{Q40} and \eqref{Q41}, we have
$$
0 \geq \lambda \int_{\mathbb{R}^3} |u|^2 dx
= -\frac{3 }{4} \int_{\mathbb{R}^3}\int_{\mathbb{R}^3} \frac{|u(y)|^2}{|x-y|} dy |u|^2 dx
+   (1 - \gamma_p) \int_{\mathbb{R}^3} |u|^p dx,
$$
which and Lemma \ref{gle4}-(2) imply
\begin{align}\label{Q42}
 (1 - \gamma_p) \int_{\mathbb{R}^3} |u|^p dx
\leq \frac{3 }{4} \int_{\mathbb{R}^3}\int_{\mathbb{R}^3} \frac{|u|^2}{|x-y|} dy |u|^2 dx
\leq C  |u|_{\frac{12}{5}}^4.
\end{align}
 By the interpolation inequality, we have
$$
  (1 - \gamma_p) |u|_p^p \leq C  |u|_{\frac{12}{5}}^4
\leq C |u|_p^{\frac{5p}{3}(1-s)} |u|_2^{\frac{10}{3}s}
\leq C   |u|_p^{\frac{5p}{3}(1-s)} a^{\frac{10}{3}s},
$$
where $s = \frac{5p-12}{5p-10}$. Due to   $p \in (3, \frac{10}{3})$, one gets $p > \frac{5p}{3}(1-s)$. Thus, the above inequality shows that
\begin{align}\label{Q50}
|u|_p^p \leq\widetilde{C} a^{\frac{2(5p-12)}{3p-8}}.
\end{align}
From \eqref{Q41}, \eqref{Q42} and \eqref{Q50}, one has
\begin{align}\label{Q57}
0=&|\nabla u|_2^2 + \frac{1}{4} \int_{\mathbb{R}^3} \int_{\mathbb{R}^3} \frac{|u|^2}{|x-y|} dy |u|^2 dx
- \int_{\mathbb{R}^3} |u|^6 dx -   \gamma_p \int_{\mathbb{R}^3} |u|^p dx\nonumber\\
\geq & |\nabla u|_2^2 - \widetilde{C}(p) a^{\frac{2(5p-12)}{3p-8}}-S^{-3} |\nabla u|_2^6:=g(|\nabla u|_2^2),
\end{align}
which gives that, there exists $0<\Upsilon_1< \frac{\sqrt{3}}{3}\mathcal{S}^{\frac{3}{2}}< \Upsilon_2<\mathcal{S}^{\frac{3}{2}}$ satisfying
\begin{align*}
|\nabla u|_2^2\leq \Upsilon_1  \quad \quad \text{or}\quad \quad  |\nabla u|_2^2\geq \Upsilon_2.
\end{align*}
If $|\nabla u|_2^2\leq \Upsilon_1 $ holds, similar to \eqref{R6}, it follows from \eqref{Q41} that
\begin{align*}
0=&|\nabla u|_2^2 + \frac{1}{4} \int_{\mathbb{R}^3}\int_{\mathbb{R}^3} \frac{|u(y)|^2}{|x-y|} dy |u|^2 dx- \int_{\mathbb{R}^3} |u|^6 dx -   \gamma_p \int_{\mathbb{R}^3} |u|^p dx \\
\geq& |\nabla u|_2^2 - C_4(p)a^{\frac{2(p-3)}{p-2}}|\nabla u|_2^{\frac{2(p-3)}{p-2}}|\nabla u|_2^2- \mathcal{S}^{-3} |\nabla u|_2^6\\
\geq& \left(1- C_4(p)a^{\frac{2(p-3)}{p-2}}  \Upsilon_1^{\frac{p-3}{p-2}}-\Upsilon_1^2 \mathcal{S}^{-3} \right)|\nabla u|_2^2\\
\geq& \left(1- C_4(p)a^{\frac{2(p-3)}{p-2}} \big(\frac{\sqrt{3}}{3}\big)^{\frac{p-3}{p-2}}\mathcal{S}^{\frac{3(p-3)}{2(p-2)}}-
\frac{1}{3} \right)|\nabla u|_2^2\\
=& \left(\frac{2}{3}- C_4(p)a^{\frac{2(p-3)}{p-2}} \big(\frac{\sqrt{3}}{3}\big)^{\frac{p-3}{p-2}}\mathcal{S}^{\frac{3(p-3)}{2(p-2)}}\right)|\nabla u|_2^2\\
>&0
\end{align*}
when $a<  3^{\frac{1}{4}}(\frac{2}{3C_4(p)})^{\frac{p-2}{2(p-3)}}\mathcal{S}^{-\frac{3}{4}}$, where
\begin{align*}
    C_4(p)=\frac{p-2}{2} \left( \frac{2(4-p)}{\overline{C}} \right)^{\frac{4-p}{p-2}}\Big(\frac{1}{\gamma_p}\Big)^\frac{2}{p-2}.
    \end{align*}
 This gives a contradiction. If $|\nabla u|_2^2\geq \Upsilon_2$,
since $\frac{1}{3}\mathcal{S}^{\frac{3}{2}}-c_pa^\eta\geq m_a= \lim_{n\rightarrow \infty} \mathcal{J}(u_n)=\lim_{n\rightarrow \infty} \big(\mathcal{J}(u_n)- \frac{1}{6}P(u_n)\big)\geq \mathcal{J}(u)- \frac{1}{6}P(u)$,
by \eqref{Q50} and $P(u)=0$, we obtain
\begin{align*}
\frac{1}{3}\mathcal{S}^{\frac{3}{2}}-c_pa^\eta\geq m_a
&\geq\mathcal{J}(u)-\frac16P(u)\\
&=\frac13|\nabla u|_2^2+\frac{5}{24}\int_{\mathbb R^3}\int_{\mathbb R^3}
\frac{|u(y)|^2|u(x)|^2}{|x-y|}\,dy\,dx
-\frac{6-p\gamma_p}{6p}|u|_p^p\\
&\geq \frac13|\nabla u|_2^2-Ca^{\frac{2(5p-12)}{3p-8}}
\end{align*}
which shows that
\begin{align}\label{Q56}
|\nabla u|_2^2\leq \mathcal S^{\frac32}-3c_pa^\eta
+Ca^{\frac{2(5p-12)}{3p-8}}:=\Theta.
\end{align}
Since $\eta<\frac{2(5p-12)}{3p-8}$, it follows that
$$
\Theta=\mathcal S^{\frac32}-3c_pa^\eta+o(a^\eta)
\qquad\text{as }a\to0^+.
$$
Consequently, recalling the definition of $g$ in \eqref{Q57}, there is a $\overline{a_*}>0$ small enough such that for every $a\in(0, \overline{a_*})$,
\begin{align*}
g(\Theta)&=\Theta-\mathcal S^{-3}\Theta^3-\widetilde C(p)a^{\frac{2(5p-12)}{3p-8}}\\
&=6c_pa^\eta+o(a^\eta)-\widetilde C(p)a^{\frac{2(5p-12)}{3p-8}}>0.
\end{align*}
This implies that $\Theta<\Upsilon_2$. Then, it follows from \eqref{Q56} that
$\Upsilon_2\leq|\nabla u|_2^2\leq\Theta<\Upsilon_2$, which is a contradiction.
Thus, we infer that $\lambda>0$.

 Next,  we let $v_n := u_n - u$, then $v_n \rightharpoonup 0$ in $\mathcal{H}_r$, $v_n \to 0$ in $L^t(\mathbb{R}^3)$ with $t \in (2,6)$ and $v_n(x) \to 0$ a.e. on $\mathbb{R}^3$. Thus,
$$
|\nabla u_n|_2^2 = |\nabla v_n|_2^2 + |\nabla u|_2^2 + o_n(1).
$$
Using the Brezis-Lieb Lemma \cite{1996W}, we have, as $n \to \infty$,
\begin{align}\label{Q43}
|u_n|_s^s = |v_n|_s^s + |u|_s^s + o_n(1),
\end{align}
where $s \in (2, 6]$. Define
$$
\Lambda(u) := \int_{\mathbb{R}^3} \int_{\mathbb{R}^3} \frac{|u(y)|^2}{|x-y|} dy|u|^2 dx.
$$
Then, we obtain from \cite[Lemma 2.2]{2008-JMAA-ZHAO} that, as $n \to \infty$,
$$
\Lambda(u_n) = \Lambda(v_n) + \Lambda(u) + o_n(1).
$$
Using the above facts, we deduce
\begin{align}\label{Q44}
\mathcal{J}(u_n) = \mathcal{J}(v_n) + \mathcal{J}(u) + o_n(1)
\end{align}
and
\begin{align}\label{Q45}
P(u_n) = P(v_n) + P(u) + o_n(1).
\end{align}
We distinguish the two cases: either (i) $\limsup_{n \to \infty} |v_n|_6^6 \to 0$ or (ii) $\limsup_{n \to \infty} |v_n|_6^6 \to l > 0$. Assuming that (ii) holds, we deduce from \eqref{Q45}, $P(u) = 0$ and the Sobolev inequality that
$$
|\nabla v_n|_2^2 = |v_n|_6^6 + o_n(1) \geq \mathcal{S}^{\frac{3}{2}} + o_n(1).
$$
Then, it follows that
$$
\mathcal{J}(v_n) \ge \frac{1}{3} \mathcal{S}^{\frac{3}{2}} + o_n(1).
$$
By \eqref{Q44},
\begin{align}\label{Q46}
\mathcal{J}(u_n) \ge \mathcal{J}(u) + \frac{1}{3} \mathcal{S}^{\frac{3}{2}} + o_n(1).
\end{align}
If $\mathcal{J}(u)\geq 0$, we get a contradiction with $m_a \leq  \frac{1}{3} \mathcal{S}^{\frac{3}{2}}-c_pa^\eta$.
Now, we claim that $\mathcal{J}(u)\geq 0$ if $a\in(0, \overline{a})$.  In fact,
 using \eqref{Q41}, and similar to  \eqref{k5}, we have
 \begin{align*}
  \mathcal{J}(u)&=\mathcal{J}(u) - \frac{1}{6} P(u)\nonumber\\
 &=\frac{1}{3} \int_{\mathbb{R}^3} |\nabla u|^2 dx + \frac{5}{24} \int_{\mathbb{R}^3} \int_{\mathbb{R}^3}  \frac{|u(y)|^2}{|x-y|} dy |u|^2 dx +  \left( \frac{\gamma_p}{6} - \frac{1}{p} \right) \int_{\mathbb{R}^3} |u|^p dx\nonumber\\
 &\geq \frac{1}{3}|\nabla u|_2^2                                                                                                                 +\frac{5}{24} \widehat{C} \frac{|u|_p^{\frac{2p}{10-3p}}}{ |\nabla u|_2^{\frac{2(3p-8)}{10-3p}} |u|_2^{\frac{8(p-3)}{10-3p}}}  -\frac{6-p\gamma_p}{6p} |u|_p^p\nonumber\\
 &\geq \left(\frac{1}{3} - C_3(p) a^{\frac{8(p-3)}{3p-8}}\right)|\nabla u|_2^2,
\end{align*}
which implies that $  \mathcal{J}(u)\geq 0$ when $a< (3 C_3(p))^{\frac{8-3p}{8(p-3)}}$.   Hence, (i) holds. Then we get from \eqref{Q43} that $|u_n|_6^6 \to |u|_6^6$ as $n \to \infty$.
Multiplying \eqref{Q29} and \eqref{Q37} by $u_n - u$ and integrating in $\mathbb{R}^3$, one gets
$$
\langle \mathcal{J}'(u_n) - \mathcal{J}'(u), u_n - u \rangle
+ \lambda_n \int_{\mathbb{R}^3} u_n (u_n - u) dx
- \lambda \int_{\mathbb{R}^3} u (u_n - u) dx
= o_n(1).
$$
By that $\lambda_n \to \lambda$, $|u_n|_6^6 \to |u|_6^6$ as $n \to \infty$ and $\mathcal{H}_r \hookrightarrow L^t(\mathbb{R}^3)$ with $t \in (2,6)$ is compact , we deduce
$$
|\nabla(u_n - u)|_2^2 + \lambda |u_n - u|_2^2 = o_n(1),
$$
which shows that $\|u_n - u\| \to 0$ as $n \to \infty$ because of $\lambda > 0$. Thus, we infer that $(\lambda, u) \in \mathbb{R}^+ \times \mathcal{H}_r$ is a solution for Eq. \eqref{L1}, where $u_n \to u$ in $\mathcal{H}_r$ and $\mathcal{J}(u) = m_a$.
\end{proof}

The following Lemma  can be found in   Br\'{e}zis and Kato   \cite[Lemma 2.1]{1979-JMPA}, which is   crucial for the $L^\infty$  estimate for the solution of Eq. \eqref{L1}.

\begin{lemma} \label{Lem4.2}
Let $\Omega \subset  \mathbb{R}^3$ and $ k(x)\in  L^{\frac{3}{2}}(\Omega )$ be a nonnegative function. Then, for every
$\varepsilon> 0$, there exists a constant $C(\varepsilon, k) > 0$ such that
$$
\int_{\Omega}k(x) u^2dx\leq \varepsilon \int_{\Omega}|\nabla u|^2dx+ C(\varepsilon, k) \int_{\Omega}| u|^2dx
$$
for all $u\in H^1(\Omega)$.
\end{lemma}

\begin{lemma}\label{Lem4.3}
Assume that  $p\in(3, \frac{10}{3} )$.
  Let $   u\in S_{a, r}$ be a nonnegative real-valued solution for Eq. \eqref{L1} and the associated Lagrange multiplier $\lambda>0$, Then there exists $C>0$ such that  $|u|_\infty <C $.
\end{lemma}

\begin{proof} Assume that $0\leq u\in S_{a, r}$ be a real-valued solution for Eq. \eqref{L1}, and the associated Lagrange multiplier $\lambda>0$, we show that $|u|_\infty <C $ for some $C>0$.
Let   $A_{m}=\{x\in\mathbb{R}^{3}:|u(x)|\leq m\}$, $B_{m}=\mathbb{R}^{3}\backslash{A_{m}}$, where  $m\in \mathbb{N}$.
For any $q>0$, define
\begin{equation*} u_m=
\begin{cases}
  |u|^{2q+1}, \quad\quad &x\in A_{m}\\
 m^{2q}u, \quad\quad&x\in B_{m}
\end{cases}
\end{equation*}
and
\begin{equation*} w_m=
\begin{cases}
|u|^{q+1},\quad\quad &x\in A_{m}\\
 m^q u,\quad\quad &x\in B_{m}.\\
\end{cases}
\end{equation*}
By simple calculation, one has
\begin{equation*} \nabla u_m=
\begin{cases}
 (2q+1)|u|^{2q}\nabla u, \quad\quad &x\in A_{m}\\
 m^{2q}\nabla u, \quad\quad&x\in B_{m}\\
\end{cases}
\end{equation*}
and
\begin{equation*} \nabla w_m=
\begin{cases}
 (q+1)|u|^{q}\nabla u, \ \ \ \ \ \ \ &x\in A_{m}\\
m^q\nabla u, \ \ \ \ \ \ \ \ \ \ \ \ \ \ &x\in B_{m}.\\
\end{cases}
\end{equation*}
Notice that $u_{m}$, $w_{m}\in \mathcal{H}$ and $w^{2}_{m}=uu_{m}\leq|u|^{2(q+1)}$.
Taking $u_{m}$ as a test function in Eq. \eqref{L1},   we get
$$
\int_{\mathbb{R}^3} \nabla u \cdot \nabla u_m+ \lambda u u_mdx +\int_{\mathbb{R}^3} \int_{\mathbb{R}^3} \frac{|u(y)|^2} {|x-y|}dyu u_{m} dx=\int_{\mathbb{R}^3} |u|^{4} u u_m+|u|^{p-2} u u_mdx,
$$
which and $\lambda>0$ imply that
\begin{align}\label{h59}
\int_{A_m} (2q+1) |u|^{2q} |\nabla u|^2dx+\int_{B_m}  m^{2q} |\nabla u|^2dx \leq \int_{\mathbb{R}^3}  |u|^{4} u u_m+|u|^{p-2} u u_mdx.
\end{align}
Since
\begin{align}\label{h60}
  \int_{\mathbb{R}^3} |\nabla w_m|^2dx=(q+1)^2 \int_{A_m}   |u|^{2q} |\nabla u|^2dx+ \int_{B_m}  m^{2q} |\nabla u|^2dx,
\end{align}
then combining \eqref{h59} with \eqref{h60}, by Lemma \ref{Lem4.2}, one has
\begin{align}
\frac{2q+1}{(q+1)^2}  \int_{\mathbb{R}^3} |\nabla w_m|^2dx&=(2q+1)\int_{A_m}   |u|^{2q} |\nabla u|^2dx+\frac{2q+1}{(q+1)^2}  \int_{B_m}  m^{2q} |\nabla u|^2dx\nonumber\\
&\leq \int_{\mathbb{R}^3}|u|^{4} u u_m+|u|^{p-2} u u_mdx\nonumber\\
&\leq \int_{\mathbb{R}^3}( |u|^{4}+|u|^{p-2} ) w_m^2dx\nonumber\\
&\leq C\int_{\mathbb{R}^3}(1+ |u|^{4} ) w_m^2dx\label{h11}\\
&\leq \varepsilon \int_{\mathbb{R}^3}   |\nabla w_m|^2dx+ C(\varepsilon,u) \int_{\mathbb{R}^3}w_m^2dx. \label{h8}
\end{align}
Furthermore, since $ w_m=|u|^{q+1}$ in $A_m$, it follows that 
\begin{align}
\Big(\int_{A_m} |  u|^{6(q+1)}dx\Big)^{ \frac{1}{3}}
= \Big(\int_{A_m} |  w_m|^{6}dx\Big)^{ \frac{1}{3}}
\leq& \Big(\int_{\mathbb{R}^3} |  w_m|^{6}dx\Big)^{ \frac{1}{3}}\nonumber\\
\leq &\mathcal{S}^{-1} \int_{\mathbb{R}^3} |\nabla w_m|^2dx\label{h9}
\end{align}
Taking $q=2$ and $\varepsilon=\frac{2}{9}$ in \eqref{h8}, and by the fact that $w_m^2\leq |u|^{2(q+1)}$ on $\mathbb{R}^3$, one has
\begin{align}\label{h10}
\int_{\mathbb{R}^3} |\nabla w_m|^2dx\leq C \int_{\mathbb{R}^3}w_m^2dx\leq C \int_{\mathbb{R}^3}|u|^{6}dx.
\end{align}
Let $m\rightarrow\infty$ in \eqref{h9}, and taking $q=2$, using \eqref{h10}  we infer that
$$
\Big(\int_{\mathbb{R}^3} |  u|^{18}dx\Big)^{ \frac{1}{3}}\leq C \int_{\mathbb{R}^3}|u|^{6}dx<+\infty,
$$
that is, $u\in L^{18}(\mathbb{R}^3)$.  Using Sobolev interpolation inequality, we get
\begin{align}\label{h12}
u\in L^{8}(\mathbb{R}^3).
\end{align}
Using \eqref{h11},  the H\"{o}lder inequality  and the Young's inequality, one has
\begin{align}\label{h13}
\frac{2q+1}{(q+1)^2}  \int_{\mathbb{R}^3} |\nabla w_m|^2dx
&\leq C\int_{\mathbb{R}^3}(1+ |u|^{4} ) w_m^2dx \nonumber \\
&\leq C\int_{\mathbb{R}^3}w_m^2dx+ C|u|_8^4 \left(\int_{\mathbb{R}^3}  w_m^4dx\right)^{\frac{1}{2}} \nonumber\\
&\leq C\int_{\mathbb{R}^3}w_m^2dx+ C|u|_8^4 \left(\int_{\mathbb{R}^3}  w_m^2dx\right)^{\frac{1}{4}} \left(\int_{\mathbb{R}^3}  w_m^6dx\right)^{\frac{1}{4}}  \nonumber\\
&\leq C\int_{\mathbb{R}^3}w_m^2dx+ C\mathcal{S}^{-\frac{3}{4}}|u|_8^4 \left(\int_{\mathbb{R}^3}  w_m^2dx\right)^{\frac{1}{4}} \left( \int_{\mathbb{R}^3} |\nabla w_m|^2dx \right)^{\frac{3}{4}} \nonumber\\
 &\leq C\int_{\mathbb{R}^3}w_m^2dx+ \delta  \int_{\mathbb{R}^3} |\nabla w_m|^2dx+
 C \delta^{-3}\mathcal{S}^{-3}|u|_8^{16} \int_{\mathbb{R}^3}  w_m^2dx,
\end{align}
where $\delta= \frac{2q+1}{4(q+1)^2} $.  Combining \eqref{h12} with \eqref{h13},    one infers that
\begin{align}\label{h14}
 \int_{\mathbb{R}^3} |\nabla w_m|^2dx\leq C(q+1)^4 \int_{\mathbb{R}^3}  w_m^2dx,
\end{align}
where $C$ is independent of $q$ and $m$. Let $m\rightarrow\infty$ in \eqref{h9},   by the fact that $w_m^2\leq |u|^{2(q+1)}$ on $\mathbb{R}^3$ and \eqref{h14}, we have
\begin{align}\label{h61}
|u|_{ {6(q+1) }}\leq C^{ \frac{1}{2(q+1)}} (q+1)^{\frac{2}{q+1}}  |u|_{ {2(q+1) } }.
\end{align}
Setting  $ \beta:=3$ and $q+1=\beta^j $ for $j\in \mathbb{N}$, then \eqref{h61} changes
$$
|u|_{6 \beta^j } \leq C^{ \frac{1}{2 \beta^j}} (\beta^j)^{\frac{2}{\beta^j} }  |u|_{6\beta^{j-1} }.
$$
Further, we proceed the $j$ times iterations that
$$
|u|_{6 \beta^j } \leq C^{\sum_{i=1}^j \frac{1}{2 \beta^i}}  \beta^{\sum_{i=1}^j\frac{2i}{\beta^i}}  |u|_{6}.
$$
Due to $\lim_{j\rightarrow\infty}C^{\sum_{i=1}^j \frac{1}{2 \beta^i}}<\infty $ and $\lim_{j\rightarrow\infty}\beta^{\sum_{i=1}^j\frac{2i}{\beta^i}} <\infty$, there exists $C>0$,   such that
\begin{align}\label{h97}
\limsup_{j\rightarrow\infty}|u|_{6 \beta^j }  \leq C.
\end{align}
We now claim that
\begin{align}\label{h98}
\lim_{j\rightarrow\infty}|u|_{6 \beta^j }=|u|_{\infty}.
\end{align}
Indeed, if $|u|_{\infty}=+\infty$, for any $n\in \mathbb{N}$, there is $\Omega_n\subset \mathbb{R}^3$ with $|\Omega_n|\in (0, +\infty)$ such that
$|u(x)|>n$ for all $x\in \Omega_n$. Hence,
$
|u|_{6 \beta^j}\geq n |\Omega_n|^{\frac{1}{6 \beta^j}}
$
implies that
$
\liminf_{j\rightarrow\infty}|u|_{6\beta^j}\geq n.
$
By the arbitrariness of $n$, we have
$
\liminf_{j\rightarrow\infty}|u|_{6\beta^j}=+\infty,
$
which is contradicts with \eqref{h97}. Thus, $|u|_{\infty}<+\infty$. Let  $j\geq1$, then
$$
\Big(\int_{\mathbb{R}^3}|u|^{6 \beta^j}dx\Big)^{\frac{1}{6 \beta^j}}\leq \Big(|u|_{\infty}^{6 \beta^j-6} \int_{\mathbb{R}^3}|u|^{6}dx\Big)^{\frac{1}{6 \beta^j}}=|u|_{\infty}^{1-\frac{1}{ \beta^j}}
|u|_{6}^{\frac{1}{\beta^j}},
$$
which implies that
$$
\limsup_{j\rightarrow\infty}|u|_{6 \beta^j}\leq |u|_{\infty}.
$$
Besides, for any $\varepsilon>0$, taking $\widetilde{\Omega}_{\varepsilon}:=\{x\in \mathbb{R}^3: |u|_{\infty}- \varepsilon\leq |u(x)|\leq  |u|_{\infty}\}$,  then $|\widetilde{\Omega}_{\varepsilon}|\in(0, +\infty)$ and
$$
\Big(\int_{\mathbb{R}^3}|u|^{6 \beta^j}dx\Big)^{\frac{1}{6\beta^j}}\geq (|u|_{\infty}-\varepsilon) |\widetilde{\Omega}_{\varepsilon}|^{\frac{1}{6 \beta^j}},
$$
which gives that
$$
\liminf_{j\rightarrow\infty}|u|_{6 \beta^j}\geq |u|_{\infty}-\varepsilon.
$$
Hence, we get \eqref{h98} by the arbitrariness of $\varepsilon$. In view of   \eqref{h97} and \eqref{h98}, one has
$
|u|_{\infty}\leq C.
$
Thus, we complete the proof.
\end{proof}

\hspace{-1cm} {\bf{Proof of Theorem \ref{K-TH1}:}}

Let
\begin{align} \label{h91}
 U_a:=\big\{ u\in\mathcal{P}_{a}:\ \mathcal{J}(u)=m_{a}=\inf_{\mathcal{P}_{a} }\mathcal{J}  \big\}.
\end{align}
From Lemmas \ref{R-lem1}-\ref{R-lem6}, one obtains that $U_a\neq\emptyset$.
Moreover, it is easy to see  from  Lemma \ref{Lem4.1}  that $  U_a= Z_{a,r}$.
We claim that
\begin{align}\label{h90}
u\in Z_{a,r} \ \Rightarrow \ |u|\in Z_{a,r}, \ \ |\nabla |u||_2=|\nabla u|_2.
\end{align}
Indeed, observe that $|\nabla |u||_2\leq|\nabla u|_2 $, then $\mathcal{J}(|u|)\leq \mathcal{J}(u)$ and $P(|u|)\leq P(u)=0$. By Lemma \ref{Q-Lem2.3}, there exists $s_{|u|}\leq 0$ satisfying $s_{|u|}\star |u|\in \mathcal{P}_{a}$. Hence, we conclude
\begin{align*}
 m_{a}&\leq \mathcal{J}(s_{|u|}\ast |u|) \\
 &=\frac{e^{2s_{|u|}}}{2}\int_{\mathbb{R}^3}|\nabla |u||^2dx+\frac{ e^{s_{|u|}}}{4} \int_{\mathbb{R}^3} \int_{\mathbb{R}^3} \frac{|u(y)|^2} {|x-y|}dy|u|^{2}dx
 -\frac{e^{6s_{|u|}}}{6}\int_{\mathbb{R}^3}|u|^{6}dx
 -\frac{ e^{p\gamma_ps_{|u|}}}{p}
 \int_{\mathbb{R}^3}|u|^{p}dx\\
 &\leq\frac{e^{2s_{|u|}}}{2}\int_{\mathbb{R}^3}|\nabla u|^2dx+\frac{ e^{s_{|u|}}}{4} \int_{\mathbb{R}^3} \int_{\mathbb{R}^3} \frac{|u(y)|^2} {|x-y|}dy|u|^{2}dx
 -\frac{e^{6s_{|u|}}}{6}\int_{\mathbb{R}^3}|u|^{6}dx
 -\frac{e^{p\gamma_ps_{|u|}}}{p}
 \int_{\mathbb{R}^3}|u|^{p}dx\\
&= \mathcal{J}(s_{|u|}\ast u)\\
&\leq  \mathcal{J}(  u)\\
&= m_{a},
\end{align*}
which shows that $|\nabla |u||_2=|\nabla u|_2 $. Thus one has $\mathcal{J}(  |u|)=\mathcal{J}(  u)=  m_{a} $ and $ P(|u|)=P(u)=0$. This gives that $|u|\in Z_{a,r}$.   Thus,  $|u|$ is a non-negative real-valued solution to Eq. \eqref{L1} for some $\lambda \in \mathbb{R}^+$ by \eqref{h85}.
Let $w:=|u|$,
  we claim that $w>0$ on $\mathbb{R}^3$.
Indeed,  from Lemma \ref{Lem4.3}, we know that
\begin{align*}
-\Delta w= -\lambda w-\left(|x|^{-1}\ast|w|^2\right)w + |w|^{2^*-2}w+ |w|^{p-2}w\in L_{loc}^q(\mathbb{R}^3,  \mathbb{R})  \quad\quad \text{for any}\ q\geq 2,
\end{align*}
which and  Cald\'{e}ron-Zygmund inequality\cite{2001book} imply $w\in W_{loc}^{2, q}(\mathbb{R}^3, \mathbb{R} )$.
 Furthermore, it follows from the Sobolev embedding theorem that $w\in C_{loc}^{1, \alpha}(\mathbb{R}^3, \mathbb{R})$ for some exponent $0 < \alpha< 1$. Applying the strong maximum principle, we further derive that $w > 0$ on $\mathbb{R}^3$. Hence,  there exists a couple $(\lambda, w)\in \mathbb{R}^+\times \mathcal{H}$ that solves equation \eqref{L1}, where $w$ is a positive real-valued function and satisfies $\mathcal{J}(w)= m_{a}$. Therefore, we complete the  proof of  Theorem \ref{K-TH1}.  $\hfill\Box$

\section{Proof of Theorems  \ref{TH4} and \ref{TH5}}\label{sec4}

 Now, we show that  the phenomenon of finite-time blow-up  occurs  for the solution of the  initial-value problem \eqref{kk1} with initial datum
  under the assumptions  of   Theorem \ref{TH4}.\\

\hspace{-1cm} {\bf {Proof of Theorem \ref{TH4}:}}    Assume that  $p\in(3,  \frac{10}{3})$ and   $\mathcal{J}(u_0)<\inf_{u\in \mathcal{P}_{a}}\mathcal{J}(u)$ with $u_0 \in S_{a, r}$.
Let $ \phi(t,x)$   be a solution of  the  initial-value problem \eqref{kk1}
with initial datum $u_0$ on $(-T_{min}, T_{max})$ satisfying
\begin{align}\label{h53}
|u_0|_2=| \phi|_2\quad\quad \text{and} \quad \quad \mathcal{J}(u_0)=\mathcal{J}( \phi).
\end{align}
  By the fact that  $|x|u_0 \in L^2(\mathbb{R}^3 ,\mathbb{C})$  and \cite[Proposition 6.5.1]{2003-C}, one gets
\begin{align} \label{h55}
H(t):=\int_{\mathbb{R}^3}|x|^2|\phi(t,x)|^2dx<+\infty  \ \ \ \ \ \text{for all}\ t\in(-T_{min}, T_{max} ).
  \end{align}
Moreover, the function $H\in C^2(-T_{min}, T_{max} )$ and the following  Virial identity holds:
\begin{align}\label{h52}
H'(t)= -4y(t)  \quad \quad \text{and} \quad \quad H''(t)=-4 y'(t)= 8 P(\phi),
  \end{align}
  where
  \begin{align}
  y(t)&= -Im  \int_{\mathbb{R}^3}\overline{\phi} (x\cdot \nabla\phi )dx;\label{h69}\\
  y'(t)
  &= -2\int_{\mathbb{R}^3}|\nabla \phi|^2dx-\frac{1}{2}\int_{\mathbb{R}^3} \int_{\mathbb{R}^3} \frac{|\phi(t,y)|^2} {|x-y|}dy|\phi(t,x)|^{2}dx+2 \int_{\mathbb{R}^3}|\phi|^{6}dx+2\gamma_p \int_{\mathbb{R}^3}|\phi|^{p}dx. \label{h70}
  \end{align}
 Now, we claim that, for any $u\in S_{a, r}$, the function $\varphi_u(s):=\psi_u(\ln s)$ with $s>0$ has a   unique global maximum point $\widehat{s}_u=e^{s_u}$ and $\varphi_u(s)$ is strictly decreasing and concave on $(\widehat{s}_u, +\infty)$.
 In fact,   from \eqref{k8} we have
  \begin{align*}
\varphi_u(s):=\psi_u(\ln s)
&=\frac{s^2}{2}\int_{\mathbb{R}^3}|\nabla u|^2dx+\frac{s}{4}  \int_{\mathbb{R}^3} \int_{\mathbb{R}^3} \frac{|u(y)|^2} {|x-y|}dy|u(x)|^{2}dx-\frac{s^{6}}{6} \int_{\mathbb{R}^3}|u|^{6}dx
-\frac{s^{p\gamma_p }}{p}\int_{\mathbb{R}^3}|u|^{p}dx,
\end{align*}
and $s \varphi_u'(s)=P(\ln s\star u)$ with $s>0$ due to \eqref{j-Q30}.  Thus, we infer that, for some $\widehat{s}_0>0$,     $\varphi_u'(\widehat{s}_0)=0$ $\Leftrightarrow$ $P(\widetilde{s}_0\star u)=0$,  where $\widetilde{s}_0=\ln \widehat{s}_0$.
Similar to Lemma \ref{Q-Lem2.3},  we obtain that, for any $u\in S_{a, r}$, $ \varphi_u'(s)=0$ has a unique solution (global
maximum point), denoted by $\widehat{s}_u$. Hence, $P(\widetilde{s}_u\star u)=0$ with $\widetilde{s}_u=\ln \widehat{s}_u$. By uniqueness, $s_u= \widetilde{s}_u$. Thus, the function $\varphi_u(s)$ with $s>0$ has a   unique global maximum point $\widehat{s}_u=e^{s_u}$. In addition, similar to Lemma \ref{Q-Lem2.3}, by simple calculation,   $\varphi_u(s)$ is strictly decreasing and concave on $(\widehat{s}_u, +\infty)$.
 According to the assumption $s_{u_0}<0$,   one obtains $ \widehat{s}_{u_0}\in(0,1)$. We claim that
\begin{align}\label{h54}
\text{if}\ u \in S_{a, r} \ \text{and}\ \widehat{s}_u\in(0,1), \ \text{then}\ P(u)\leq \mathcal{J}(u)-\inf_{\mathcal{P}_{a}}\mathcal{J}.
\end{align}
In fact, since $ \widehat{s}_u\in(0,1)$ and $\varphi_u(s)$ is strictly decreasing and concave on $(\widehat{s}_u, +\infty)$, one shows that $P(u)<0$ and
\begin{align*}
\mathcal{J}(u)=\psi_u(0) = \varphi_u(1)\geq  \varphi_u(\widehat{s}_u )-\varphi_u'(1)(\widehat{s}_u-1 )
=\mathcal{J}(s_u \ast u) -|P(u)|(1-\widehat{s}_u )
\geq \inf_{\mathcal{P}_{a}}\mathcal{J} +P(u),
\end{align*}
which yields \eqref{h54} holds.
Now, let us consider the solution $\phi$  for the  initial-value problem \eqref{kk1} with initial datum $u_0$. Since by assumption $s_{u_0} < 0$, and the map $u \mapsto s_u$ is continuous, we deduce that  $s_{ \phi(t)} < 0$   for every $|t|<\overline{t}$ with $\overline{t}>0$ small enough. Then $\widehat{s}_{ \phi(t)}\in (0,1)$ for $|t|<\overline{t}$. By \eqref{h54}  and recalling the assumption     $\mathcal{J}(u_0)< \inf_{\mathcal{P}_{a}}\mathcal{J} $, it follows  from \eqref{h53} that
$$
P(\phi(t))\leq \mathcal{J}(\phi(t) )-\inf_{\mathcal{P}_{a}}\mathcal{J}=\mathcal{J}(u_0 )-\inf_{\mathcal{P}_{a}}\mathcal{J}:=-\eta<0,
$$
for every such $|t|< \overline{t}$. Next, we show that
\begin{align}\label{h56}
  P(\phi(t) )\leq-\eta \quad\quad \text{ for any }\ t\in (-T_{min}, T_{max}).
\end{align}
 Assume that there is $t_0\in (-T_{min}, T_{max})$ satisfying $ P(\phi(t_0) )=0$. Then, by \eqref{h53} one has $\phi(t_0)\in S_{a, r}$  and
$$
 \mathcal{J}(\phi(t_0) )\geq m_{a}=\inf_{\mathcal{P}_{a}}\mathcal{J}> \mathcal{J}(u_0)=  \mathcal{J}(\phi(t_0) ),
$$
which is a contradiction.  This shows that    $ P(\phi(t_0) )\neq0$ for any  $t_0\in (-T_{min}, T_{max})$. Then, since $P(\phi(t))<0$ for every   $|t|< \overline{t}$, we obtain $ P(\phi(t) )<0$ for any  $t\in (-T_{min}, T_{max})$ (if at some $t\in  (-T_{min}, T_{max})$, $P(\phi(t))>0$. By continuity, we have $P(\phi(\widetilde{t}))=0$ for some $\widetilde{t}\in  (-T_{min}, T_{max})$, a contradiction). Since $ P(\phi(t) )<0$ for any  $t\in (-T_{min}, T_{max})$, by Lemma \ref{Q-Lem2.3} we have $s_{\phi(t) }<0$. That is, $\widehat{ s}_{\phi(t) }\in (0,1)$. Therefore, \eqref{h54} holds and   the above arguments yield
\begin{align*}
  P(\phi(t) )\leq-\eta \quad\quad \text{ for any }  t\in (-T_{min}, T_{max}).
\end{align*}
Thus, we deduce from \eqref{h55}, \eqref{h52}, \eqref{h56} and Taylor's formula that
\begin{align*}
0\leq H(t)= H(0)+H'(0) t+\frac{1}{2}H''( \xi) t^2\leq H(0)+H'(0) t-4 \eta t^2
\end{align*}
for any $t\in (-T_{min}, T_{max})$ and some $\xi\in (-T_{min}, T_{max})$,  which implies that
 $T_{max}$ has   an upper bound since the right hand side becomes negative for $t$ large.
 That is,
the    solution $\phi$  of the initial-value problem  \eqref{kk1}
with initial datum $u_0$   blows-up in finite time. \\

In the following, we are going to prove   the strong instability of the standing wave for problem \eqref{kk1} by  using   Theorem \ref{TH4}.

 \hspace{-0.9cm}{\bf {Proof of Theorem \ref{TH5}:}}
According to \eqref{h90}, we know that if $\widehat{u}$ is a ground state radial solution for \eqref{L1}, then $|\widehat{u}|$ is also a ground state radial solution to \eqref{L1}. Thus, without loss of generality, assume that $\widehat{u}>0$. Similar to \eqref{h85}, the associated Lagrange multiplier $\widehat{\lambda}>0$.
 For any $p\in (3, \frac{10}{3})$, let    $\phi_\varrho(x, t)$ be the solution to   the  initial-value problem \eqref{kk1} with initial datum $\widehat{u}_\varrho$, where $\widehat{u}_\varrho=\varrho\star\widehat{u}$ with $\varrho>0$.  We have $\widehat{u}_\varrho \rightarrow \widehat{u} $ in $\mathcal{H}$ as $\varrho\rightarrow 0^+$, and hence it is sufficient to prove
that $\phi_\varrho(x, t)$  blows-up in finite time. Let $s_{\widehat{u}_\varrho}\in \mathbb{R}$ be defined in Lemma \ref{Q-Lem2.3}. Clearly $s_{\widehat{u}_\varrho} = -\varrho < 0$, and
$$
 \mathcal{J}(\widehat{u}_\varrho)= \mathcal{J}(\varrho\star \widehat{u})<\mathcal{J}(  \widehat{u})=\inf_{\mathcal{P}_{a}}\mathcal{J}.
$$
 Since $\widehat{u}\in \mathcal{H}_r$, by radial Lemma (see \cite{1983-BL}), we have
\begin{align}\label{h68}
\widehat{u}(x) \to 0\ \ \ \ \  \mathrm{as }~|x|\to \infty.
\end{align}
Thus, for any $p\in(3, \frac{10}{3})$, by \eqref{h68}  and the fact that $\widehat{\lambda}>0$, there exists    $\widehat{R}>0$ large enough such that
 \begin{align*}
-\Delta \widehat{u}=\big[-(|x|^{-1}\ast |\widehat{u}|^2)-\widehat{\lambda}+|\widehat{u}|^{4}+|\widehat{u}|^{p-2}\big ]\widehat{u}\leq -\frac{\widehat{\lambda}}{2}\widehat{u}     \quad \quad \text{for all}\ |x|>\widehat{R}.
 \end{align*}
Let $\xi(x)=M_1e^{\big(-\sqrt{\frac{\widehat{\lambda}}{2}}|x|\big)}$, where $M_1$ satisfies
 \begin{align*}
M_1e^{-\sqrt{\frac{\widehat{\lambda}}{2}}\widehat{R}}\geq \widehat{u}(x)\quad\quad \text{for}\ |x|=\widehat{R}.
 \end{align*}
By simple calculation,  we get $ \Delta \xi \leq\frac{\widehat{\lambda}}{2} \xi$ for all $x\neq0$. Set $\zeta=\xi-\widehat{u}$, one has
 \begin{align*}
 \begin{cases}
-\Delta \zeta+\frac{\widehat{\lambda}}{2} \zeta \geq 0,\quad &|x|>\widehat{R},\\
\zeta(x)\geq0,\quad &|x|=\widehat{R},\\
\lim_{|x|\rightarrow\infty}\zeta(x)=0.
 \end{cases}
 \end{align*}
Hence, it follows from maximum principle \cite{2001book} that $\zeta\geq0$. That is,
$$
\widehat{u}(x)\leq M_1e^{\big(-\sqrt{\frac{\widehat{\lambda}}{2}}|x|\big)}\quad\quad \text{ for all} \ |x|\geq \widehat{R}.
$$
Moreover, since $\widehat{u}$ is continuous function, there exists $M_2>0$ such that
$$
\widehat{u}(x) e^{\big(\sqrt{\frac{\widehat{\lambda}}{2}}|x|\big)} \leq M_2 \quad\quad \text{ for all} \ |x|\leq\widehat{R}.
$$
Therefore, choosing $M=\max\{ M_1, M_2\}$, we have
\begin{align}\label{k6}
\widehat{u}(x)\leq Me^{\big(-\sqrt{\frac{\widehat{\lambda}}{2}}|x|\big)}\quad\quad \text{ for all} \ x\in \mathbb{R}^3.
\end{align}
Then we deduce from \eqref{k6} that
$$
\int_{\mathbb{R}^3} |x|^2|\widehat{u}_{\varrho}|^2dx=e^{-2 \varrho} \int_{\mathbb{R}^3} |x|^2|\widehat{u}|^2dx<+\infty,
$$
which shows that
$ |x|\widehat{u}_\varrho\in L^2(\mathbb{R}^3, \mathbb{C})$.
Hence, it follows from Theorem \ref{TH4}  that the solution $\phi_\varrho(x, t)$  blows-up in finite time. Thus, by  Definition \ref{de1.2}, the
standing wave $e^{i\widehat{\lambda} t}\widehat{u} $ is strongly unstable. We complete the proof of Theorem \ref{TH5}.$\hfill\Box$\\

\noindent{\bf Statements and Declarations }

 {\bf   Competing Interests:}~The authors declare that they have no known competing financial interests or personal
relationships that could have appeared to influence the work reported in this paper.\\

\noindent{\bf{Data Availability Statements}}

Data sharing not applicable to this article as no datasets were generated or analysed during the current study.\\

\noindent{\bf{Acknowledgments}}

We sincerely appreciate the anonymous reviewers for their valuable comments and suggestions, which have significantly improved the quality of our manuscript.
This work is supported by  National Natural Science Foundation of China (No.~12371120 and No.~12661024); and the National Natural Science Foundation of China under the Tianyuan Fund for Mathematics (No.~12526512); and the Science and Technology Research Program of Chongqing Municipal Education Commission (No. KJQN202601320).


\begin{thebibliography}{50}


\setlength{\itemsep}{-2.5pt}

\bibitem{2010AIHP-ADP}
  \newblock  Azzollini, A.,  d'Avenia, P.,  Pomponio, A.:
   \newblock  On the Schr\"{o}dinger-Maxwell equations under the effect of a general nonlinear
term.
\newblock  Ann. Inst. H. Poincar\'{e} Anal. Non Lin\'{e}aire  \textbf{27}, 779--791 (2010)

\bibitem{2011ZAMP-BS}
\newblock  Bellazzini, J.,  Siciliano,  G.:
\newblock Stable standing waves for a class of nonlinear Schr\"{o}dinger-Poisson
equations.
\newblock Z. Angew. Math. Phys. \textbf{62},  267--280 (2011).

\bibitem{2011JFA-BS}
\newblock  Bellazzini, J.,  Siciliano,  G.:
\newblock  Scaling properties of functionals and existence of constrained minimizers.
\newblock J. Funct. Anal. \textbf{261}, 2486--2507 (2011).


\bibitem{2013PLMS}
\newblock  Bellazzini, J.,  Jeanjean, L.,   Luo, T.:
\newblock Existence and instability of standing waves with prescribed norm for a class of
Schr\"{o}dinger-Poisson equations.
\newblock Proc. Lond. Math. Soc. \textbf{107},  303--339 (2013)


\bibitem{2002RMP-BF}
 \newblock  Benci, V.,   Fortunato, D.:
  \newblock Solitary waves of the nonlinear Klein-Gordon equation coupled with Maxwell equations.
   \newblock Rev. Math. Phys.  \textbf{14}, 409--420   (2002)



\bibitem{1983-BL}
\newblock  Berestycki, H.,   Lions, P.L.:
 \newblock Nonlinear scalar field equations I, existence of a ground state.
 \newblock Arch. Ration. Mech. Anal.   \textbf{82},  313--345 (1983)



\bibitem{1979-JMPA}
\newblock  Brezis, H.,   Kato, T.:
  \newblock Remarks on the Schr\"{o}dinger operator with singularly complex potentials.
  \newblock J. Math. Pures Appl. \textbf{58},  137--151 (1979)

\bibitem{1983-CPAM-Brezis}
\newblock  Br\'{e}zis, H.,  Nirenberg,  L.:
\newblock Positive solutions of nonlinear elliptic equations involving critical Sobolev exponents.
\newblock  Commun. Pure Appl. Math. \textbf{36},  437--477 (1983)




\bibitem{2003-C}
\newblock  Cazenave, T.: Semilinear Schr\"{o}dinger Equations, Courant Lecture Notes in Mathematics, vol. 10, New York University, Courant Institute of Mathematical Sciences/American Mathematical Society, New York/Providence, RI, (2003)




\bibitem{2013-3M}
\newblock Catto I., Dolbeault J., S\'{a}nchez \'{O}.,  Soler  J.:
 \newblock Existence of steady states for the Maxwell-Schr\"{o}dinger-Poisson system: exploring the applicability of the concentration-compactness principle.
 \newblock Math. Models Methods Appl. Sci. \textbf{23}, 1915--1938 (2013)








\bibitem{2020JMAA-CT}
\newblock  Chen, S.T.,  Tang, X.H., Yuan, S.:
 \newblock Normalized solutions for Schr\"{o}dinger-Poisson equations with general nonlinearities.
\newblock J. Math. Anal. Appl. \textbf{481},   123447 (2020)



\bibitem{2022.9}
\newblock  Chen, S.T.,  Tang, X.H.:
\newblock   Another look at Schr\"{o}dinger equations with prescribed mass.
 \newblock J. Differential Equations \textbf{386},  435--479 (2024)




\bibitem{2004-ANS-D}
\newblock  D'Aprile, T.,   Mugnai,  D.:
\newblock Non-existence results for the coupled Klein-Gordon-Maxwell equations.
\newblock  Adv. Nonlinear Stud. \textbf{4},  307--322 (2004)






\bibitem{1993book}
 \newblock  Ghoussoub, N.:
 \newblock Duality and Perturbation Methods in Critical Point Theory, Cambridge Tracts in Mathematics, vol.107, Cambridge University Press, Cambridge, (1993).
%
%
\bibitem{2001book}
\newblock  Gilbarg, D.,    Trudinger, N.S.:  Elliptic Partial Differential Equations of Second Order. Reprint of the 1998 edition. Classics in Mathematics. Springer-Verlag, Berlin (2001)


\bibitem{2004-ANS-Ha}
  \newblock  Hajaiej, H.,   Stuart, C.A.:
 \newblock  On the variational approach to the stability of
standing waves for the nonlinear Schr\"{o}dinger equation.
  \newblock  Adv. Nonlinear Stud.   \textbf{ 4},  469--501 (2004)



  \bibitem{1997jean}
\newblock  Jeanjean, L.:
  \newblock Existence of solutions with prescribed norm for semilinear elliptic equations.
   \newblock  Nonlinear Anal.  \textbf{28}, 1633--1659  (1997)


  \bibitem{2021-JDE-JL}
  \newblock  Jeanjean, L.,   Le T.T.:
  \newblock Multiple normalized solutions for a Sobolev critical Schr\"{o}dinger-Poisson-Slater equation.
  \newblock J. Differential Equations \textbf{303}, 277--325  (2021)



\bibitem{2013ZAMP-JL}
\newblock  Jeanjean, L.,  Luo, T.:
\newblock Sharp nonexistence results of prescribed $L^2$-norm solutions for some class of Schr\"{o}dinger-Poisson and quasi-linear equations,
\newblock  Z. Angew. Math. Phys. \textbf{64}, 937--954 (2013)


\bibitem{2023-JGA-KLT}
\newblock  Kang, J.C., Liu, X.Q., Tang, C.L.:
\newblock  Ground state sign-changing solutions for critical Schr\"{o}dinger-Poisson system with steep potential well.
\newblock  J. Geom. Anal.   \textbf{33}, 59 (2023)


\bibitem{kang}
  \newblock Kang, J.C., Li, Y.Y., Tang, C.L.:
  \newblock Prescribed mass standing waves for Schr\"{o}dinger-Maxwell equations with combined nonlinearities. Discrete Contin. Dyn. Syst. \textbf{45},  3296--3344 (2025)



\bibitem{2025-JDE-Kang}
\newblock Kang, J.C., Li, Y.Y., Tang, C.L.:
\newblock Multiple normalized solutions for Schr\"{o}dinger-Maxwell equation with Sobolev critical exponent and mixed nonlinearities.
\newblock  J. Differential Equations \textbf{443},  113564 (2025)

\bibitem{2013APDE}
\newblock  Lenzmann, E.,  Lewin, M.:
\newblock Dynamical ionization bounds for atoms.
\newblock Anal. PDE  \textbf{6}, 1183--1211 (2013)




\bibitem{2017ARMA}
\newblock   Killip, R.,   Oh, T.,  Pocovnicu, O.,   Visan,  M.:
 \newblock Solitons and scattering for the cubic-quintic nonlinear
Schr\"{o}dinger equation on $\mathbb{R}^3$.
\newblock Arch. Rational Mech. Anal.  \textbf{225},   469--548 (2017)

\bibitem{2026-Li-JGA}
\newblock Li, M.R., He, X.M.:
\newblock Multiple Normalized Solutions for the Quasilinear Sobolev Critical Schr\"{o}dinger-Poisson System.
\newblock J. Geom. Anal. \textbf{36}, 231 (2026).


\bibitem{1984CMP-L}
 \newblock   Lions, P.L.:
 \newblock  Solutions of Hartree-Fock equations for Coulomb systems.
   \newblock  Comm. Math. Phys. \textbf{109}, 33--97  (1984)


\bibitem{2016-AMPA-LWZ}
  \newblock Liu, Z.L., Wang, Z.Q., Zhang, J.J.:
   \newblock Infinitely many sign-changing solutions for the nonlinear
Schr\"{o}dinger-Poisson system.
 \newblock  Ann. Mat. Pura Appl. \textbf{195}, 775--794 (2016)

\bibitem{2026-Liu-JGA}
\newblock Liu, W.,  He, X.M.,  Meng,  Y.X.:
 \newblock Prescribed $L^2$-norm solutions for a non-autonomous Schr\"{o}dinger-Poisson system with critical growth.
\newblock  J. Geom. Anal. \textbf{36},  23 (2026)




\bibitem{2014JMAA-L}
\newblock  Luo, T.J.:
\newblock Multiplicity of normalized solutions for a class of nonlinear Schr\"{o}dinger-Poisson-Slater equations.
\newblock J. Math. Anal. Appl. \textbf{416}, 195--204  (2014)



\bibitem{2013CMP}
\newblock  Miao, C.,   Xu,  G.,     Zhao, L.:
\newblock The dynamics of the 3D radial NLS with the combined terms.
\newblock Commun. Math. Phys.  \textbf{318},   767--808 (2013)



\bibitem{2025-CV-P}
\newblock  Peng, X.,   Rizzi, M.:
\newblock  Normalized solutions of mass supercritical Schr\"{o}dinger-Poisson equation with potential.
\newblock  Calc. Var. Partial Differential Equations  \textbf{64}, 152 (2025)




\bibitem{2006JFA}
\newblock   Ruiz,  D.:
\newblock     The Schr\"{o}dinger-Poisson equation under the effect of a nonlinear local term.
\newblock     J. Funct. Anal.  \textbf{237},  655--674 (2006)


\bibitem{2004JSP}
\newblock  Sanchez,  O.,  Soler, J.:
\newblock  Long-time dynamics of the Schr\"{o}dinger-Poisson-Slater system.
 \newblock  J. Statist. Phys. \textbf{114}, 179--204  (2004)





\bibitem{2023IJM}
 \newblock  Siciliano, G.,   Silva, K.:
   \newblock   On the structure of the Nehari set associated to a Schr\"{o}dinger-Poisson system with prescribed mass: old and new results.
   \newblock   Israel J. Math. \textbf{258}, 403--451 (2023)

\bibitem{2020Soave}
 \newblock  Soave, N.:
  \newblock Normalized ground state for the NLS equations with combined nonlinearities.
  \newblock   J. Differential Equations  \textbf{269},  6941--6987 (2020)


\bibitem{2021-JFA-Soave}
\newblock Soave, N.:
  \newblock Normalized ground states for the NLS equation with combined nonlinearities: The Sobolev critical case.
  \newblock  J. Funct. Anal.  \textbf{279}, 108610  (2020)


\bibitem{2007CPDETao}
\newblock  Tao, T.,   Visan, M., Zhang,  X.:
\newblock The nonlinear Schr\"{o}dinger equation with combined power type nonlinearities.
\newblock  Commun. Partial Differ. Equ.  \textbf{32}, 1281--1343 (2007)




\bibitem{2022JFAWei}
\newblock  Wei, J.C.,   Wu, Y.Z.:
\newblock Normalized solutions for Schr\"{o}dinger equations with
critical Sobolev exponent and mixed nonlinearities.
\newblock J. Funct. Anal.  \textbf{283},  109574 (2022)


\bibitem{1983W}
\newblock  Weinstein, M.I.:
 \newblock Nonlinear Schr\"{o}dinger equations and sharp interpolation estimates.
 \newblock Comm. Math. Phys.  \textbf{87},   567--576  (1983)

\bibitem{1996W}
 \newblock   Willem, M.:
  \newblock  Minimax Theorems, vol. 24, Birkh\"{o}user, Boston, Mass., (1996)

\bibitem{2018ZAMP-YL}
\newblock  Ye, H.Y.,    Luo, T.J.:
\newblock On the mass concentration of $L^2$-constrained minimizers for a class of Schr\"{o}dinger-Poisson equations.
\newblock Z. Angew. Math. Phys.   \textbf{69},  66 (2018)

\bibitem{2017JMAA-ZZ}
\newblock Zeng,  X.Y.,  Zhang, L.:
\newblock Normalized solutions for Schr\"{o}dinger-Poisson-Slater equations with unbounded potentials.
\newblock J. Math. Anal. Appl. \textbf{452}, 47--61  (2017)


\bibitem{2008-JMAA-ZHAO}
\newblock  Zhao, L.G., Zhao, F.K.:
\newblock  On the existence of solutions for the Schr\"{o}dinger-Poisson equations.
\newblock   J. Math. Anal. Appl.  \textbf{346},    155--169 (2008)








\end{thebibliography}

\end{document}